\documentclass[11pt]{amsart}
\usepackage[T1]{fontenc}
\usepackage[latin1]{inputenc}
\usepackage{a4wide}
\usepackage{setspace}
\usepackage{xypic}
\usepackage{amssymb}
\usepackage{amsthm}
\usepackage[french]{babel}
\usepackage{latexsym}

\usepackage{srcltx}
\date{}

\newtheorem{thm}{Th\'eor\`eme}[section]
\newtheorem*{thm*}{Th\'eor\`eme}

\newtheorem{defn}[thm]{D\'efinition}
\newtheorem{rem}[thm]{Remarque}
\newtheorem{prop}[thm]{Proposition} 
\newtheorem{lem}[thm]{Lemme} 
\newtheorem{cor}[thm]{Corollaire} 
\newcommand{\Coprod}{\displaystyle\coprod}
\newcommand{\Prod}{\displaystyle\prod}
\newcommand{\fonc}[5]{
 \begin{array}{cccc}
 #1: & #2 & \longrightarrow & #3\\
     & #4 & \longmapsto & #5
 \end{array}
}

\begin{document}

\title{Utilisation d'une cohomologie étale équivariante en topologie arithmétique.}
\author{Baptiste Morin}
\maketitle
\renewcommand{\abstractname}{Abstract}


\begin{abstract}
A. Sikora has shown in \cite{Sikora} results which confirm the
analogy between number fields and 3-manifolds. However, he has given
proofs of his results which are very different in the arithmetic and
in the topological case. In this paper, we show how to provide a
unified approach to the results in the two cases. For this we
introduce an equivariant cohomology which satisfies a localization
theorem. In particular, we obtain a satisfactory explanation for the
coincidences between Sikora's formulas which leads us to clarify and
to extend the dictionary of arithmetic topology.
\end{abstract}

\section{Introduction.}

Le calcul de la cohomologie étale du spectre $X$ de l'anneau
d'entiers d'un corps de nombres $L$ laisse imaginer une analogie
entre les corps de nombres et les variétés réelles de dimension
trois. En effet, lorsque $L$ est totalement imaginaire, les groupes
de cohomologie étale de $X$ à coefficients dans un faisceau abélien
arbitraire sont nuls après la dimension trois (cf \cite{Deninger}
4.6). De plus, le théorème de dualité d'Artin-Verdier, analogue
arithmétique de la dualité de Poincaré, fait à nouveau apparaître la
dimension trois comme dimension maximale. Par ailleurs, le spectre
d'un corps fini est, du point de vue de la topologie étale, un objet
de dimension un dont le groupe fondamental est le complété profini
de $\mathbb{Z}$. Pour cette raison, les points fermés de $X$ sont
vus comme des noeuds dans une variété de dimension trois. Une
extension galoisienne $L/K$ de groupe $G$ correspond à un revêtement
galoisien $M\rightarrow M/G$ de variétés compactes de dimension
trois. Dans cette situation, le groupe des classes $Cl(L)$ et le
quotient libre $U_L/\mu_L$  du groupe des unités de $L$
correspondent respectivement, en tant que modules galoisiens, à la
partie de torsion $H_{tor}(M)$ et au quotient libre $H_{free}(M)$ du
premier groupe d'homologie singulière de $M$ à coefficients entiers.
Cette analogie est connue sous le nom de topologie arithmétique et
nous renvoyons à \cite{Reznikovpreprint} pour un dictionnaire plus
complet.

Lorsque le groupe de Galois est le groupe cyclique d'ordre premier
$C_p$ et qu'il opère par automorphismes préservant l'orientation, le
lieu de branchement du revêtement $M\rightarrow M/C_p$ est constitué
d'un nombre fini de noeuds ramifiés, analogues topologiques des
places finies ramifiées dans l'extension $L/K$. Cette situation est
envisagée dans \cite{Sikora}. A partir de la structure galoisienne
des groupes $H_{tor}(M)$ et $H_{free}(M)$ (respectivement $Cl(L)$ et
$U_L/\mu_L$), A.~Sikora y donne un encadrement du nombre $s$ de
noeuds (respectivement de places finies) ramifiés. Il obtient des
résultats en accord quasi-parfait avec le dictionnaire de la topologie
arithmétique. Cependant, ses preuves sont basées sur des méthodes
très différentes, puisqu'il utilise une cohomologie équivariante
dans le contexte géométrique alors qu'il fait appel à à la théorie
du corps de classes en arithmétique.

Le but de ce travail est de fournir des démonstrations analogues à
ces résultats analogues afin de comprendre la coïncidence à priori
très surprenante de ces résultats. Alors, ces résultats, leurs
hypothèses et leurs démonstrations mettent en interaction la
quasi-totalité des éléments du dictionnaire de la topologie
arithmétique, et permettent ainsi de tester la validité de ces
correspondances. Afin  d'approfondir ce dictionnaire, la pertinence
de certains de ses éléments et les contradictions offertes par
d'autres sont alors mises en évidence.

La cohomologie équivariante modifiée utilisée dans \cite{Sikora}
pour traiter le cas des variétés topologiques satisfait
essentiellement deux propriétés. Elle est d'une part l'aboutissement
d'une suite spectrale dont le terme initial est de la forme
$\widehat
H^p(G;H^q(M;\mathbb{Z}))\Longrightarrow\widehat{H}^{p+q}_G(M;\mathbb{Z}).$
D'autre part, cette cohomologie équivariante satisfait un théorème
de localisation, et par conséquent, ne "voit" que le lieu de
ramification. Dans la section 2, nous définissons l'analogue en
cohomologie étale de cette théorie et la suite spectrale $\widehat
H^p(G;H^q(X_{et};F))\Longrightarrow\widehat{H}^{p+q}_G(X_{et};F)$
qui y aboutit. Les groupes $\widehat{H}^{q}_G(X_{et};F)$ ne sont
définis que lorsque le groupe d'automorphismes $G$ est fini, et ne
peuvent être intéressants que si les groupes de cohomologie de $F$
sont nuls à partir d'une certaine dimension $n$. Dans ce cas, les
groupes $\widehat{H}^{q}_G(X_{et};F)$ s'identifient à partir de la
dimension $n+1$ aux groupes de cohomologie mixte $H^q(X_{et};G;F)$
(cf \ref{comparaison}), définis comme les invariants cohomologiques
du topos des $G$-faisceaux d'ensembles sur $X_{et}$ (cf
\cite{classedechern} 2).


Nous démontrons dans la section 3 le théorème de localisation énoncé
ci-dessous.
\begin{thm*}
Soit $X$ un schéma noethérien sur lequel un groupe fini $G$ opère
fidèlement et de manière admissible. On note $Z$ le lieu de
ramification du revêtement $X\rightarrow X/G$ et
$\phi:\widetilde{Z}\rightarrow X$ la limite projective des
voisinages étales de $Z$. Soit de plus $F$ un $G$-faisceau adapté
(\ref{adapté}) sur $X$. Alors le morphisme canonique
$\widehat{H}_{G}^{*}(X_{et};F) \longrightarrow
\widehat{H}_{G}^{*}(\widetilde{Z}_{et};\phi^*F)$ est un
isomorphisme. Si de plus $F$ est de torsion et si $Z$ est contenu
dans un ouvert affine, on a l'isomorphisme
$\widehat{H}_{G}^{*}(X_{et};F) \simeq
\widehat{H}_{G}^{*}(Z_{et};i^*F)$, où $i:Z\rightarrow X$ désigne
l'immersion fermée canonique.
\end{thm*}
Ce théorème permet par exemple de calculer les groupes
$\widehat{H}^q_G(X_{et};\mathbb{G}_m)$, lorsque $X$ est le spectre
de l'anneau d'entiers d'un corps de nombres $L$. On trouve en
particulier lorsque $G$ est abélien, l'isomorphisme
$\widehat{H}^0_G(X_{et};\mathbb{G}_m)\simeq\Prod I_{\mathfrak{q}},$
où le produit est pris sur l'ensemble des places $\mathfrak{q}$ de
$L^G$ et pour lesquelles $I_{\mathfrak{q}}$ désigne le sous-groupe
d'inertie dans $G$. La suite spectrale établit donc un lien entre la
ramification dans l'extension $L/L^G$ et la structure galoisienne
des groupes $Cl(L)$ et $U_L$. La première moitié de la section 4 est
consacrée à l'étude de cette suite spectrale et à certaines de ses
conséquences. On obtient en particulier une démonstration des
résultats de A. Sikora en théorie des nombres à l'aide de ces
méthodes.

Nous exprimons dans la section 4.4 la relation de dualité donnant
les isomorphismes
$$\widehat{H}_G^n(X_{et};\mathbb{G}_m)\simeq\widehat{H}_G^{2-n}(X_{et};\mathbb{Z})^D,$$
d'ailleurs compatibles à la dualité induite par celle
d'Artin-Verdier sur les termes initiaux des suites spectrales
respectives. L'utilisation du $G$-faisceau $\mathbb{Z}$ permet alors
de donner des preuves satisfaisantes du point de vue de la topologie
arithmétique.


En effet, la preuve des résultats de \cite{Sikora} de nature
topologique et arithmétique respectivement s'articule comme suit. Il
s'agit d'encadrer le nombre~$s$ de noeuds (respectivement de places
finies) ramifiés dans un revêtement $X\rightarrow X/C_p$ de variétés
de dimension trois (respectivement de spectres d'anneaux d'entiers
de corps de nombres) galoisien de groupe cyclique d'ordre
premier~$p$. Le théorème de localisation fournit les isomorphismes
$$\widehat{H}_{C_p}^{n}(X;\mathbb{Z}) \simeq
 \widehat{H}_{C_p}^{n}(Z;\mathbb{Z})\simeq\mathbb{F}_p^s,$$ où $Z$
désigne le lieu de ramification qui est constitué de $s$ noeuds
(respectivement de $s$ places finies). La suite spectrale
$$\widehat
H^i(C_p;H^j(X;\mathbb{Z}))\Longrightarrow\widehat{H}^{i+j}_{C_p}(X;\mathbb{Z})$$
donne ainsi des approximations successives du module
$\mathbb{F}_p^s$ à partir des groupes $\widehat
H^i(C_p;H^j(X;\mathbb{Z}))$. La dualité de Poincaré (respectivement
d'Artin-Verdier) permet alors d'obtenir un encadrement du nombre $s$
en fonction de la dimension sur $\mathbb{F}_p$ des espaces
$\widehat{H}^0(C_p;H_{tor}(X))$ et $\widehat{H}^1(C_p;H_{free}(X))$
(respectivement $\widehat{H}^0(C_p;Cl(L))$ et
$\widehat{H}^1(C_p;U_L/\mu)$).

L'intérêt de ces groupes de cohomologie étale équivariante modifiée
en topologie arithmétique est qu'ils sont les stricts analogues des
mêmes groupes définis dans le contexte topologique, ce qui n'est pas
le cas des premiers groupes de cohomologie équivariante non
modifiée. Ceci vient renforcer l'analogie des preuves esquissées
ci-dessus.

Cependant, l'utilisation de la topologie étale d'Artin-Verdier est
nécessaire pour traiter le cas des extensions de corps de nombres
non totalement imaginaires.

Pour terminer, nous tirons les conclusions de ce travail dans une
cinquième section. Nous prenons alors clairement parti pour la
première des deux versions (sensiblement différentes) du
dictionnaire de la topologie arithmétique proposées par A. Reznikov
dans \cite{Reznikovpreprint} et \cite{Reznikov}. On est ainsi amené
à montrer comment la cohomologie de S. Lichtenbaum (conjecturale à
ce jour (cf \cite{MatFlach})), associée à la topologie Weil-étale,
fournit naturellement un analogue arithmétique au groupe
$H_1(M;\mathbb{Z})$, dont le sous-groupe de torsion et le quotient
libre maximal s'identifient, en tant que modules galoisiens, à
$Cl(L)$ et $U_L/\mu_L$ respectivement. Les calculs de S. Lichtenbaum
confirment donc clairement l'intuition des fondateurs de la
topologie arithmétique. Cependant, il semble que l'analogue
arithmétique du groupe fondamental $\pi_1(M)$ ne puisse être le
groupe de Galois $G_L^{nr}$ de l'extension maximale non ramifiée de
$L$.

Ces considérations suggèrent que les résultats de \cite{Sikora} sont
deux manifestations d'un même phénomène, l'arithmétique apparaissant
ici comme un "cas particulier" du cadre topologique.

\

\emph{Remerciements:} Je tiens à remercier Boas Erez pour m'avoir
initié à la topologie arithmétique et pour m'avoir proposé cette
problématique. Je suis également très reconnaissant envers Lorenzo
Ramero. Ses lectures attentives et ses nombreuses remarques ont été
très utiles tout au long de ce travail. Je tiens aussi à remercier
Matthias Flach, pour m'avoir transmis ses résultats non publiés sur
la cohomologie du groupe de Weil.

\section{Cohomologie étale équivariante modifiée.}
On fixe un schéma $X$ sur lequel un groupe fini $G$ opère (à gauche)
par automorphismes. Dans la suite, on appelle faisceau sur $X$ un
faisceau de groupes abéliens pour la topologie étale et on note
$Ab(X)$ la catégorie des faisceaux sur $X$.
\subsection{La catégorie des $G$-faisceaux.}

Nous rappelons quelques propriétés de cette catégorie indispensables
pour la suite et pour lesquelles nous n'avons pas de référence.

\begin{defn}
Soient $X$ un schéma muni d'une action d'un groupe fini $G$ et $F$
un faisceau sur $X$. On appelle \emph{$G$-linéarisation} de $F$ la
donnée d'une famille $(\varphi_{\sigma})_{\sigma\in G}$ de
morphismes de faisceaux $\varphi_{\sigma}:\sigma_*F\rightarrow F$,
telle que les conditions suivantes soient satisfaites:
\begin{itemize}
\item $\varphi_1=Id$.
\item $\varphi_{\tau\sigma}=\varphi_{\tau}\circ {\tau}_{*}(\varphi_{\sigma})$ (condition
de cocycle).
\end{itemize}
\end{defn}

Un \emph{$G$-faisceau} est un faisceau muni d'une $G$-linéarisation.
Un \emph{morphisme de $G$-faisceaux} $\alpha:F\rightarrow L$ est un
morphisme de faisceaux $\alpha:F\rightarrow L$ tel que les
diagrammes
 \[ \xymatrix{
\sigma_{*}F \ar[d]_{\sigma_{*}(\alpha)} \ar[r]^{\varphi_{F;\sigma}} & F \ar[d]_{\alpha}   \\
\sigma_{*}L \ar[r]^{\varphi_{L;\sigma}}        & L } \]  soient
commutatifs, les structures de $G$-faisceaux sur $F$ et $L$ étant
définies par les $G$-linéarisations $(\varphi_{F;\sigma})_{\sigma\in
G}$ et $(\varphi_{L;\sigma})_{\sigma\in G}$. Les faisceaux usuels en
topologie étale (par exemple $\mathbb{G}_{m}$ ou les faisceaux
constants) sont naturellement munis de $G$-linéarisations (cf
\cite{tohoku}5.1).

On note $Ab(G;X)$ la catégorie des $G$-faisceaux sur $X$ et de leurs
morphismes.

\begin{rem}
Etant donné deux $G$-faisceaux $F$ et $L$ sur $X$, on fait
opérer $G$ sur le groupe $Hom_{Ab(X)}(F;L)$ de la manière suivante.
Si $\alpha:F\rightarrow L$ est un morphisme de faisceaux et $\sigma$
un élément de $G$, on pose
$\sigma*\alpha:=\varphi_{L;\sigma}\circ\sigma_*(\alpha)\circ\varphi_{F;\sigma}^{-1}$.
Alors $Hom_{Ab(G;X)}(F;L)$ est le sous-groupe des invariants de
$Hom_{Ab(X)}(F;L)$ sous l'action de $G$.
\end{rem}
On montre alors facilement que la catégorie $Ab(G;X)$ est additive,
puisque $Ab(X)$ l'est. De plus, si $\alpha:F\rightarrow L$ est un
morphisme de $G$-faisceaux, la $G$-linéarisation de $F$ en induit
une sur $Ker(\alpha)$, et celle définie sur $L$ en induit une sur
$Coker(\alpha)$. D'autre part, un isomorphisme de $G$-faisceaux est
un morphisme de $G$-faisceaux qui est un isomorphisme en tant que
morphisme de faisceaux. Donc si $\alpha$ est un morphisme dans
$Ab(G;X)$, le morphisme $\overline{\alpha}:Coim(\alpha)\rightarrow
Im(\alpha)$ est un isomorphisme de $G$-faisceaux. Autrement dit,
$Ab(G;X)$ est une catégorie abélienne.

\begin{rem}
Soit $\phi:X\rightarrow Y$ un morphisme de schémas sur lesquels le
groupe $G$ opère. On dit que $X$ est muni d'une action de $G$
compatible à celle définie sur $Y$ lorsque $\phi$ commute à l'action
de $G$. Dans ces conditions, $\phi_*$ est un foncteur de la
catégorie des $G$-faisceaux sur $X$ dans celle des $G$-faisceaux sur
$Y$. En effet, si $F$ est un $G$-faisceau sur $X$, la
$G$-linéarisation de $F$ est transportée sur $\phi_*F$ par le
foncteur $\phi_*$. De même, $\phi^ *$ est un foncteur exact de
$Ab(G;Y)$ dans $Ab(G;X)$. De plus, si $F$ et $L$ sont des
$G$-faisceaux sur $X$ et $Y$ respectivement, l'isomorphisme
d'adjonction
\begin{equation}\label{adjonction}
Hom_{Ab(X)}(\phi^*F;L)\longrightarrow Hom_{Ab(X)}(F;\phi_*L)
\end{equation}
commute à l'action de $G$ définie ci-dessus, car le morphisme
d'ajonction $F\rightarrow\phi_*\phi^*F$ est un morphisme de
$G$-faisceaux. L'isomorphisme (\ref{adjonction}) identifie donc les
sous-groupes des invariants $Hom_{Ab(G;X)}(\phi^*F;L)$ et
$Hom_{Ab(G;X)}(F;\phi_*L)$. Ainsi, les foncteurs $\phi^*$ et
$\phi_*$ entre les catégories de $G$-faisceaux sur $Y$ et $X$, sont
adjoints. Il suit que $\phi_*$ préserve les $G$-faisceaux injectifs.
\end{rem}

\begin{prop}\label{oubliinjectif}
La catégorie $Ab(G;X)$ possède suffisamment d'injectifs. De plus, un
$G$-faisceau injectif est aussi injectif en tant que faisceau.
\end{prop}
\begin{proof}
On vérifie facilement que le foncteur
$$\fonc{Sym}{Ab(X)}{Ab(G;X)}{L}{\sum g_{*}L=\prod g_{*}L}$$ est adjoint à gauche
du foncteur d'oubli $Ou:Ab(G;X)\rightarrow Ab(X)$. De plus, les
foncteurs $g_*$ et $\sum$ sont exacts, donc $Sym$ l'est aussi. On en
déduit que le foncteur d'oubli préserve les injectifs.

D'autre part, $Sym$ est aussi adjoint à droite du foncteur d'oubli.
Si $F$ est un $G$-faisceau, on peut choisir un faisceau injectif $I$
dans lequel se plonge $F$. Ceci induit un morphisme injectif de
$G$-faisceaux $F\rightarrow Sym(I)$. De plus, l'isomorphisme
d'adjonction
$$Hom_{Ab(G;X)}(-;Sym(I))\simeq Hom_{Ab(X)}(-;I)\circ Ou$$
montre que $Sym(I)$ est un $G$-faisceau injectif. En effet, le
foncteur d'oubli est exact et $I$ est injectif. Ainsi, tout
$G$-faisceau se plonge dans un $G$-faisceau injectif.
\end{proof}

Si $I^{*}(F)$ est une résolution injective de $F$ dans la catégorie
$Ab(G;X)$, alors $I^{*}(F)(X)$ est un complexe de
$\mathbb{Z}[G]$-modules. Il suit que les groupes de cohomologie
usuels sont des $\mathbb{Z}[G]$-modules (à gauche).

Le fait que la catégorie $Ab(G;X)$ possède suffisamment d'injectifs
peut se démontrer plus directement. La preuve suivante est inspirée
par celle donnée par Grothendieck dans \cite{tohoku}. Elle nécessite
la notion de $G$-système de points géométriques (cf \ref{syst}).

\begin{prop}
La catégorie $Ab(G;X)$ possède suffisamment d'injectifs.
\end{prop}
\begin{proof}

Soit $F$ un $G$-faisceau sur $X$. On se donne un $G$-système de
points géométriques $M$ sur $X$ en gardant les notations de la
définition précédente. Pour chaque orbite $T$ de $X'$ sous l'action
de $G$, on choisit un point $\overline{\beta}$. On a la suite de
morphismes $\xymatrix{\overline{\beta} \ar[r]^{i_{\beta}}
&T\ar[r]^{v_T}&X} $, où $v_T$ est induit par $u:X'\rightarrow X$ et
$i_{\beta}$ est l'inclusion de $\overline{\beta}$ dans $T$. Soient
$G_{\beta}$ le stabilisateur du point géométrique $\beta$ et
$F_{\beta}$ la fibre de $F$ en $\beta$. Le
$\mathbb{Z}[G_{\beta}]$-module $F_{\beta}$ se plonge dans un
$\mathbb{Z}[G_{\beta}]$-module injectif $I_{\beta}$. Par ailleurs,
le stabilisateur de $g(\beta)$ est $G_{g(\beta)}:=gG_{\beta}g^{-1}$.
On définit alors $I_{g(\beta)}$ comme le groupe abélien $I_{\beta}$
sur lequel opère $G_{g(\beta)}$ à travers l'isomorphisme
$G_{g(\beta)}\rightarrow G_{\beta}$. C'est un
$\mathbb{Z}[G_{g(\beta)}]$-module injectif dans lequel se plonge
$F_{g(\beta)}$. On définit de cette manière un foncteur
$\xymatrix{Ab(G_{\beta};\overline{\beta})
\ar[r]^{Ind_{G_{\beta}}^{G}}&Ab(G;T)}$, qui est une équivalence de
catégories et dont le foncteur quasi-inverse est $i_{\beta}^{*}$. La
suite $\xymatrix{Ab(G_{\beta};\overline{\beta})
\ar[r]^{Ind_{G_{\beta}}^{G}}&Ab(G;T) \ar[r]^{v_{T *}}&Ab(G;X)}$
montre que chaque $I_{\beta}$ définit un $G$-faisceau injectif sur
$X$, puisque ces deux foncteurs préservent les injectifs. On définit
$I$ comme le produit
$\Prod_{T}v_{T*}(Ind_{G_{\beta}}^{G}(I_{\beta}))$, qui est un
$G$-faisceau injectif dans lequel se plonge $F$.

\end{proof}

\begin{prop}
Si $\pi:X\rightarrow Y$ est un revêtement étale galoisien de groupe
$G$, la catégorie des $G$-faisceaux sur $X$ est équivalente à celle
des faisceaux sur $Y$.
\end{prop}

\begin{proof}
Si $A$ est un faisceau sur $Y$, alors $\pi^*A$ est un $G$-faisceau
sur $X$. En effet, $G$ opère trivialement sur $Y$ et $A$, et $\pi$
commute à l'action de $G$. D'autre part, si $F$ est un $G$-faisceau
sur $X$, on définit un faisceau $\pi_*^G F$ sur $Y$ en posant, pour
tout $Y$-schéma étale $U$, $$\pi_*^G F(U):=F(X\times_YU)^G.$$ Le
morphisme $X\times_YU\rightarrow U$ est un revêtement étale (non
nécessairement connexe) galoisien de groupe $G$ quel que soit $U$
étale sur $Y$. Alors, si $A$ est un faisceau sur $Y$ et $U$ un
$Y$-schéma étale, le groupe $A(U)$ s'identifie à $A(X\times_YU)^G$
(cf \cite{Milne} II.1.4). On déduit facilement de ceci que les
foncteurs $\pi^*$ et $\pi_*^G$ sont quasi-inverses l'un de l'autre.
\end{proof}

\subsection{Cohomologie équivariante.}
On fixe une résolution complète $W_*$ pour le groupe fini $G$ (cf
\cite {Br} VI.3 ou \cite{cart} XII.3). C'est un complexe de
$\mathbb{Z}[G]$-modules (à gauche) libres de type fini et tel que la
cohomologie du complexe $Hom_{\mathbb{Z}[G]}(W_*;M)$ donne les
groupes de cohomologie modifiés $\widehat{H}^{*}(G;M)$, pour tout
$\mathbb{Z}[G]$-module $M$.

\begin{defn} \label{cohomologie équivariante}
Soient $F$ un $G$-faisceau et $0\rightarrow F\rightarrow
I^0\rightarrow I^1\rightarrow ...$ une résolution injective de $F$
dans $Ab(G;X)$. Alors $0\rightarrow I^0(X)\rightarrow
I^1(X)\rightarrow...$ est un complexe de $\mathbb{Z}[G]$-modules. On
note $Hom_{\mathbb{Z}[G]}(W_*;I^{*}(X))$ le double complexe
d'homomorphismes et
$$Tot^n(W_*;I^{*}(X)):=\bigoplus_{i+j=n}Hom_{\mathbb{Z}[G]}(W_i;I^{j}(X))$$
le complexe (de cochaînes) total associé. La différentielle totale y
est définie comme dans \cite {Weibel} 2.7.4. On définit les
\emph{groupes de cohomologie étale équivariante modifiée de $X$ à
coefficients dans $F$} de la manière suivante:
$$\widehat{H}^{*}_G(X;F):=H^{*}(Tot^*(W_*;I^{*}(X))).$$
\end{defn}
Ces groupes sont en fait définis comme l'aboutissement de la suite
spectrale \ref{suite spectrale}. C'est ce qui justifie l'utilisation
de la somme directe (et non du produit direct) dans la définition du
complexe total. Cependant, il semble que cette définition ne donne
des résultats intéressants que lorsqu'elle est appliquée à des
$G$-faisceaux possédant une résolution flasque finie, c'est-à-dire
quand la définition précédente du complexe total $Hom$ prend du
sens.
\begin{prop}
La cohomologie du complexe $Tot(W_*;I^{*}(X))$ ne dépend pas de la
résolution injective de $F$ choisie.
\end{prop}
\begin{proof}
Soient $I^{*}$ et $J^{*}$ deux résolutions injectives de $F$ dans
$Ab(G;X)$. Il existe des morphismes de complexes $f:I^{*}\rightarrow
J^{*}$ et $h:J^{*}\rightarrow I^{*}$ qui relèvent le morphisme
$Id:F\rightarrow F$. Alors $h\circ f$ et $Id_{(I^{*})}$ relèvent
l'identité et sont homotopes d'après \cite {Weibel} 2.7.4. De même,
$f\circ h$ et $Id_{(J^{*})}$ sont homotopes. En appliquant à ceci le
foncteur des sections globales, on obtient les morphismes de
complexes de $\mathbb{Z}[G]$-modules $f_X:I^{*}(X)\rightarrow
J^{*}(X)$ et $h_X:J^{*}(X)\rightarrow I^{*}(X)$ de sorte que
$f_X\circ h_X$ et $h_X\circ f_X$ soient homotopes à l'identité du
complexe de $\mathbb{Z}[G]$-modules correspondant. D'autre part,
deux morphismes homotopes entre deux doubles complexes induisent les
mêmes morphismes sur les groupes de cohomologie des complexes totaux
associés (cf \cite {cart} 15.6.1). Les flèches $f_X$ et $h_X$
induisent donc des isomorphismes réciproques entre les groupes
$H^{*}(Tot(W_*;I^{*}(X)))$ et $H^{*}(Tot(W_*;J^{*}(X)))$.

\end{proof}

\begin{prop} \label{suite spectrale}
Le groupe gradué $\widehat{H}^{*}_G(X;F)$ est l'aboutissement d'une
suite spectrale dont le terme initial est $E_2^{p;q}(X;F)=\widehat
H^p(G;H^q(X;F))$.
\end{prop}
\begin{proof}
La première filtration du double complexe $Hom(W_*;I^{*}(X))$ (que
l'on note provisoirement $D^{**}$) est régulière (cf \cite {cart}
15.6). La première suite spectrale de ce double complexe converge
donc vers $\widehat{H}^{*}_G(X;F):=H^{*}(Tot(W_*;I^{*}(X)))$. Le
terme initial de cette suite spectrale est
$E_2^{p;q}=H_h^p(H_v^q(D^{**}))$, où $H_h$ et $H_v$ désignent la
cohomologie du double complexe $D^{**}$ relativement aux
différentielles horizontales et verticales. Le foncteur d'oubli
$Ou:Ab(G;X)\rightarrow Ab(X)$ est exact et préserve les injectifs
(cf \ref{oubliinjectif}). L'égalité
$E_2^{p;q}=\widehat{H}^p(G;H^q(X;F))$ en découle.
 \end{proof}

Cette suite spectrale ainsi que les groupes de cohomologie
équivariante sont fonctoriels en $F$ et $X$. Plus précisemment, un
morphisme de $G$-faisceaux $f:F\rightarrow S$ sur $X$ induit un
morphisme
$$f_G^{*}:\widehat{H}^{*}_G(X;F)\rightarrow\widehat{H}^{*}_G(X;S)$$
et un morphisme de suites spectrales
$$f^{**}_*:E^{**}_*(X;F)\rightarrow E^{**}_*(X;S).$$ Ces deux
morphismes sont compatibles et $f_2^{p;q}$ s'identifie au morphisme
canonique $$\widehat{H}^p(G;H^q(X;F))\rightarrow
\widehat{H}^p(G;H^q(X;F)).$$ De plus, si $$0\rightarrow
F'\rightarrow F\rightarrow F''\rightarrow 0$$ est une suite exacte
de $G$-faisceaux, on a la suite exacte longue de cohomologie
équivariante
$$...\rightarrow\widehat{H}^{n}_G(X;F')\rightarrow\widehat{H}^{n}_G(X;F)\rightarrow\widehat{H}^{n}_G(X;F'')\rightarrow\widehat{H}^{n+1}_G(X;F')\rightarrow...$$

En effet, soient $I_{F'}^*$, $I_{F}^*$ et $I_{F''}^*$ des
$G$-résolutions injectives respectivement de $F'$, $F$ et $F''$
telles que le mophisme de complexes $$0\rightarrow
I_{F'}^*\rightarrow I_{F}^*\rightarrow I_{F''}^*\rightarrow 0$$ soit
exact. En lui appliquant le foncteur des sections globales, on
trouve la suite exactes de complexes de $\mathbb{Z}[G]$-modules
$$0\rightarrow I_{F'}^*(X)\rightarrow I_{F}^*(X)\rightarrow
I_{F''}^*(X)\rightarrow 0$$ et donc la suite exacte
 $$0\rightarrow Tot^n(W_*;I_{F'}^{*}(X))\rightarrow Tot^n(W_*;I_{F}^{*}(X))\rightarrow Tot^n(W_*;I_{F''}^{*}(X))\rightarrow 0,$$ car les $W_i$ sont projectifs. La suite exacte longue de cohomologie s'en déduit.

Soient $X$ et $Y$ deux schémas sur lesquels opère le groupe $G$ et
$F$ un $G$-faisceau sur $Y$. Un morphisme $\pi:X\rightarrow Y$
compatible à l'action de $G$ induit un morphisme
$$\pi^*_G:\widehat{H}^{*}_G(Y;F)\rightarrow\widehat{H}^{*}_G(X;f^*F)$$
et de manière compatible, un morphisme de suites spectrales
$$\pi_*^{**}:E^{**}_*(Y;F)\rightarrow E^{**}_*(X;f^*F)$$ tel que
$\pi_2^{p;q}$ s'identifie à la flèche canonique
$$\widehat{H}^p(G;H^q(Y;F))\rightarrow
\widehat{H}^p(G;H^q(X;f^*F).$$ En effet, si $I^*$ et $J^*$ sont des
$G$-résolutions de $F$ et $f^*F$ respectivement, $\pi_G^*$ et
$\pi_*^{**}$ sont induits par le morphisme de complexes
$$I^*(Y)\rightarrow f^*(I^*)(X)\rightarrow J^*(X).$$

La multiplication par le cardinal de $G$ annule toute la suite
spectrale, y compris les groupes de cohomologie équivariante (cf
\cite{Swan} 1.2). Si le groupe $G$ opère trivialement sur $X$ et
$F$, la suite spectrale $E_*^{*;*}(X;F)$ est triviale, c'est à dire
$E_2^{*;*}(X;F)=E_{\infty}^{*;*}(X;F)$ (cf \cite{Swan} 1.2).

Un faisceau sur $X$ est dit  \emph{flasque} si $\check H^q(\{ U_i
\rightarrow U \}_i;F)=0$ pour tout $q \geq1$ et tout recouvrement
$\{ U_i \rightarrow U \}_i$ (pour la topologie étale) d'un
$X$-schéma étale $U$. D'après \cite {Milne} 3.2.12, $F$ est flasque
si et seulement si $H^q(U;F)=0$ pour tout $q \geq1$ et tout
$X$-schéma étale $U$. De tels faisceaux sont acycliques pour le
foncteur des sections globales (cf \cite {Milne} 3.1.8). On appelle
\emph{$G$-résolution acyclique de $F$} toute $G$-résolution de $F$
(i.e. toute résolution de $F$ dans $Ab(G;X)$) par des $G$-faisceaux
acycliques pour le foncteur des sections globales.

\begin{prop} \label {acyclique}
Toute $G$-résolution acyclique de $F$ permet de calculer les groupes
de cohomologie étale équivariante de $X$ à coefficients dans $F$.
\end{prop}
\begin{proof}
Soit $0\rightarrow F\rightarrow I^0\rightarrow I^1\rightarrow ...$
une résolution injective de $F$ dans $Ab(G;X)$ et  $0\rightarrow
F\rightarrow C^0\rightarrow C^1\rightarrow ...$ une $G$-résolution
acyclique de $F$. On a un morphisme de complexes $f:C^{*}\rightarrow
I^{*}$ qui relève l'identité d'ailleurs unique à homotopie près. En
appliquant le foncteur des sections globales, on obtient le
morphisme de complexe de $\mathbb{Z}[G]$-modules
$f_X:C^{*}(X)\rightarrow I^{*}(X)$, qui donne des isomorphismes sur
les groupes de cohomologie. Le morphisme de complexes $f_X$ induit
un morphisme de doubles complexes $Hom(W_*;C^{*}(X))\rightarrow
Hom(W_*;I^{*}(X))$, donc un morphisme au niveau des suites
spectrales (convergentes) associées, qui est un isomorphisme dès la
deuxième page. Ce qui précède définit un et un seul (car $f$ est
défini à homotopie près) isomorphisme
$H^{*}(Tot(W_*;C^{*}(X)))\simeq H^{*}(Tot(W_*;I^{*}(X)))$ (cf \cite
{cart} 15.3.2).
\end{proof}

\begin{rem}\label{comparaison}
Soit $F$ un $G$-faisceau sur $X$ dont les groupes de cohomologie
$H^q(X;F)$ sont nuls pour $q\geq n+1$. Sous cette hypothèse, on
construit dans la section suivante une $G$-résolution acyclique de
$F$ de longueur $n$, qui d'après le théorème précédent, permet
d'obtenir les groupes $\widehat{H}^q_G(X;F)$. Soit $0\rightarrow
F\rightarrow C^0\rightarrow C^1\rightarrow ...\rightarrow
C^n\rightarrow 0$ cette résolution.

La suite
$$...\rightarrow W_{i+1}\rightarrow W_i\rightarrow ...\rightarrow
W_0 \rightarrow\mathbb{Z}\rightarrow0$$ est une résolution de
$\mathbb{Z}$ par des $\mathbb{Z}[G]$-modules projectifs. Les groupes
de cohomologie du complexe total associé au double complexe
$(Hom_{\mathbb{Z}[G]}(W_i;C^j(X)))_{i,j\geq0}$ (situé dans le
premier quadrant) satisfont les conditions axiomatiques des
foncteurs dérivés droits du foncteur (composé), qui à un
$G$-faisceau $F$, associe le groupe $F(X)^G$ des sections globales
invariantes sous l'action de $G$. Il s'agit donc des groupes de
cohomologie mixte $H^*(X;G;F)$, définis comme les invariants
cohomologiques du topos des $G$-faisceaux d'ensembles sur $X$ (cf
\cite{classedechern} 2).

On observe que les deux complexes totaux définissant respectivement
les groupes $\widehat{H}^*_G(X;F)$ et $H^*(X;G;F)$ coincident en
dimension supérieure ou égale à $n$. Donc pour $q\geq n+1$, on a
l'identification
$$\widehat{H}^q_G(X;F)\simeq H^q(X;G;F).$$
\end{rem}

%

\subsection{La résolution flasque de Godement.}
Soit $F$ un faisceau sur $X$. La résolution de Godement est une
résolution de $F$ par des faisceaux flasques. Lorsque $X$ est muni
d'une action d'un groupe fini $G$ et que $F$ est un $G$-faisceau sur
$X$, cette résolution doit être définie à partir d'un système de
points géométriques stable sous l'action de $G$, afin de conserver
une strucure équivariante. Il s'agit alors d'une $G$-résolution
flasque qui permet, grâce à la proposition \ref {acyclique}, de
calculer la cohomologie étale équivariante de $X$. Nous verrons dans
la section 3 comment cette résolution peut être utilisée pour
démontrer un théorème de localisation.

\begin{defn} \label{syst}
On appelle \emph{$G$-système de points géométriques sur $X$} un
ensemble $M$ de points géométriques
$\alpha:\overline{\alpha}\rightarrow X$ satisfaisant les conditions
suivantes.

-Si $x\in X$, il existe un point géométrique de $M$ dont l'image est
$x$.

-Si $\alpha\in M$ et $g\in G$, alors $g(\alpha):=g\circ \alpha \in
M$.

-$M$ est minimal pour ces propriétés.

On pose $X':=\Coprod_{\alpha\in M}\overline{\alpha}$. Les morphismes
$(\alpha)_{\alpha\in M}$ induisent un morphisme $u:X'\rightarrow X$
compatible à l'action de $G$ sur $X$ et $X'$.
\end{defn}
On construit un $G$-système de points géométriques de la manière
suivante. Pour toute trajectoire $T$ de $G$ sur $X$, on choisit un
point $x\in T$ et un point géométrique $\alpha_T$ dont l'image est
$x$. On pose alors $X':=\Coprod g(\alpha_T)$, où la somme est
indexée sur l'ensemble des trajectoires $T$ et sur l'ensemble des
éléments $g$ du groupe $G$.

On garde les notations de la définition précédente.
\begin{defn}On définit la résolution de Godement (équivariante)
$$0\rightarrow F\rightarrow C^0(F)\rightarrow C^1(F)\rightarrow...$$de manière
récurrente comme suit.
\begin{enumerate}
\item On pose $C^{0}(F):=u_{*}u^{*}(F)$. On a un morphisme canonique injectif
$\varepsilon:F\rightarrow C^{0}(F)$.
\item On pose $Z^{1}(F):=Coker(\varepsilon)$ et $C^{1}(F):=C^{0}(Z^1(F))$. Il y a un morphisme canonique
$d^{0}:C^{0}(F)\rightarrow C^{1}(F)$.
\item On définit par récurrence $Z^{n}(F):=Coker(d^{n-2})$ et
$C^{n}(F):=C^{0}(Z^n(F))$.
 Il y a un morphisme canonique
$d^{n-1}:C^{n-1}(F)\rightarrow C^{n}(F)$.
\end{enumerate}
\end{defn}

Les faisceaux $C^n(F)$ sont flasques. De plus, $u$ est compatible à
l'action de $G$, la catégorie $Ab(G;X)$ est abélienne et les
foncteurs $u^*$ et $u_*$ sont adjoints. Il s'agit donc d'une
$G$-résolution flasque.

Si $F$ est un $G$-faisceau sur $X$ et si $\varphi:U\rightarrow X$
est un morphisme étale, on a $C^0(F)(U)=\Prod_{\varphi\circ\beta\in
M}F_{\beta},$ où le produit est pris sur l'ensemble des points
géométriques $\beta$ de $U$ tels que $\varphi\circ\beta$ soit un
élément de $M$.

Il sera utile dans la suite de disposer de $G$-résolutions flasques finies.
Dans ce but, on définit pour tout entier $n$ la résolution de
Godement de $F$ tronquée au cran $n$
$$C_{(n)}^*(F):C^0(F)\rightarrow C^1(F)\rightarrow ... \rightarrow Im(d^{n-1}:C^{n-1}(F)\rightarrow C^n(F))\rightarrow 0.$$
Si pour tout $X$-schéma étale $U$ et pour tout $q\geq n+1$, on a
$H^q(U;F)=0$, alors $C_{(n)}^*(F)$ est une résolution flasque finie
de $F$. En effet, les $C^i(F)$ sont flasques et on a
$$H^r(U;Im(d^{n-1}))=H^{n+r}(U;F)=0$$ pour tout $X$-schéma étale $U$
et pour tout $r\geq 1$. Le faisceau $Im(d^{n-1})$ est donc flasque.

\section{Le théorème de localisation.}

Nous démontrons un théorème permettant, sous certaines hypothèses,
de calculer la cohomologie étale équivariante modifiée d'un schéma
$X$ en se restreignant au sous-schéma fermé constitué des points de
$X$ dont le groupe d'inertie est non trivial. Il s'agit de
l'analogue d'un théorème de localisation pour les espaces
topologiques de dimension finie (cf \cite{Godement} V.12)
initialement démontré par Swan dans \cite{Swan}. Cette hypothèse de
finitude se traduit ici par le fait que ce nouveau théorème de
localisation ne s'applique qu'à des $G$-faisceaux "adaptés" en un
sens que nous allons définir.

 La preuve de ce théorème comporte trois étapes. On commence par montrer
que les groupes $\widehat{H}^{*}_G(X;F)$ sont nuls, dès que le
groupe $G$ opère sur $X$ sans inertie et que $F$ est adapté. Ceci
permet ensuite de se concentrer aux $G$-voisinages étales des points
ramifiés. On obtient alors le théorème en passant à la limite sur
ces derniers.

\subsection{Le cas non ramifié.}
Soit $F$ un faisceau sur $X$ possédant une résolution flasque finie,
disons de longueur $n$. Alors les groupes $H^q(U;F)$ sont nuls pour
tout $X$-schéma étale $U$ et tout $q\geq n+1$. Réciproquement, on a
remarqué dans la section précédente que cette condition permettait
de construire une résolution flasque finie. Cette observation motive
la définition suivante.

\begin{defn} \label{adapté}
Soit $X\rightarrow Y$ un revêtement étale galoisien de groupe $G$ et
$F$ un $G$-faisceau sur $X$. On dit que $F$ est \emph{adapté} si il
 existe un entier $n$ de sorte que pour tout $Y$-schéma étale $V$ et pour tout $q\geq n+1$ on ait $H^q(V;\pi_*^GF)=0$.
\end{defn}
Cette condition de finitude sera nécessaire pour utiliser le lemme
ci-dessous.
\begin{lem} \label{Ra}
Soit $M^*$ un complexe de cochaines de $\mathbb{Z}[G]$-modules tel
que $M^{n}=0$ pour $n$ assez grand et pour $n$ négatif. On suppose
de plus que $\widehat{H}^{i}(G;M^{j})=0$ pour tout $i$ et $j$. Alors
on a $H^{*}(Tot(W_*;M^*))=0$.
\end{lem}
\begin{proof}
Dans ces conditions, la seconde filtration du double complexe
$Hom_{\mathbb{Z}G}(W_*;M^*)$ est régulière (\cite{cart}II.15.6),
donc la seconde suite spectrale converge vers $H^{*}(Tot(W_*;M^*))$.
Le terme $E_{1}$ de cette suite spectrale est
$\widehat{H}^{i}(G;M^{j})=0$, d'où $H^{*}(Tot(W_*;M^*))=0$.
\end{proof}
\begin{prop} \label{net}
Soit $\pi:X\rightarrow Y$ un revêtement étale galoisien de groupe
$G$ avec $Y$ localement noethérien. Si $F$ est un $G$-faisceau
adapté sur $X$, alors $\widehat{H}^{*}_G(X;F)=0$.
\end{prop}
\begin{proof}
On note abusivement $0\rightarrow \pi_*^GF\rightarrow C^0\rightarrow
C^1\rightarrow ...\rightarrow C^n\rightarrow 0$ la résolution de
Godement de $\pi_*^GF$ tronquée au cran $n$. C'est une résolution
flasque de $\pi_*^GF$. Comme $\pi^{*}$ est exact et qu'il préserve
les faisceaux flasques, $0\rightarrow \pi^{*}\pi_*^GF\rightarrow
\pi^{*}C^0\rightarrow \pi^{*}C^1\rightarrow ...\rightarrow
\pi^{*}C^n\rightarrow 0$ est une $G$-résolution flasque de
$\pi^{*}\pi^G_*F=F$. On veut montrer que les $C^j(X)=\pi ^*C^j(X)$
sont des $\mathbb{Z}[G]$-modules cohomologiquement triviaux.

Soit $S$ un sous groupe de $G$. Le morphisme $X\rightarrow X/S$ est
un revêtement étale galoisien de groupe $S$ et $\phi :X/S\rightarrow
Y$ est un morphisme étale (cf \cite{SGA1} V.3.3). Les faisceaux
$\phi ^*C^j$ sont flasques et on a pour tout $q$ strictement positif
(cf \cite{Milne} III.2.6): $$H^q(S;C^j(X))=\check H^q(\{X\rightarrow
X/S\};\phi ^*C^j)=0,$$ où l'on considère les groupes de cohomologie
de $\check{C}ech$ relatifs au recouvrement étale $\{X\rightarrow
X/H\}$. Pour tout sous-groupe $S$ de $G$, les $C^j(X)$ sont donc
$\Gamma^S$-acycliques. On en déduit que les $C^j(X)$ sont
cohomologiquement triviaux (cf \cite{Serre} IX.thm 8).
 En appliquant le lemme \ref {Ra}, on obtient $\widehat{H}^{*}_G(X;F)=0$.

\end{proof}

\begin{rem}
Sous les hypothèses de la proposition précédente, la suite spectrale
$$H^p(G;H^q(X;F))\Rightarrow H^{p+q}(G;X;F)$$ est celle
d'Hochschild-Serre (cf \cite{Milne} III.2.20), et les groupes de
cohomologie mixte $H^*(G;X;F)$ s'identifient aux $H^*(Y;\pi_*^GF)$.
Ces groupes sont en général non nuls. Afin de ne décrire que la
ramification, et de permettre ainsi l'existence d'un théorème de
localisation, les groupes $\widehat{H}_G^*(X;F)$ doivent s'annuler
dans cette situation.
\end{rem}

On dit qu'un groupe fini $G$ opère sur un schéma $X$ de façon
\emph{admissible} si $X$ est réunion d'ouverts affines invariants
par $G$ ou encore si toute trajectoire de $G$ dans $X$ est contenue
dans un ouvert affine. Le quotient $X/G$ existe à cette condition
(cf \cite {SGA1} V.1.7). Soit $X$ un schéma sur lequel un groupe
fini $G$ opère de manière admissible et $\varphi:U\rightarrow X$ un
morphisme étale. On suppose que $U$ est muni d'une action de $G$
compatible à celle définie sur $X$, c'est à dire que $\varphi$
commute à l'action de $G$. Si de plus le morphisme $\varphi$ est
affine, l'action de $G$ sur $U$ est admissible. En effet, soit
$(X_i)_i$ un recouvrement de $X$ par des ouverts affines de $X$
stables par $G$, alors $(\varphi^{-1}(X_i))_i$ est un recouvrement
de $U$ possédant les mêmes propriétés. De même, si
$\varphi:U\rightarrow X$ est une immersion ouverte, l'action de $G$
sur $U$ est admissible. En effet, soit $T$ une trajectoire de $U$ et
$W$ un ouvert affine de $X$ la contenant. On a $T\subset U\cap W
\subset W$. Comme dans $W$, toute partie finie a un système
fondamental de voisinages ouverts affines, il existe un voisinage
ouvert affine de $T$ contenu dans $U\cap W$ donc dans $U$.

Si un groupe fini $G$ opère sur $X$ de manière admissible et sans
inertie (i.e. tous les groupes d'inertie sont triviaux), alors
$X\rightarrow X/G$ est un revêtement étale galoisien de groupe $G$
(cf \cite{Raynaud} X corollaire 1 et \cite{SGA1} V.1.8). Si de plus
$X$ est supposé localement noethérien, $X/G$ l'est aussi (cf
\cite{Raynaud} X corollaire 2).

\begin{cor} \label{NET}
Soit $X$ un schéma localement noethérien sur lequel un groupe $G$
opère de manière admissible et sans inertie. Soient de plus $F$ un
$G$-faisceau adapté sur $X$ et $U$ un $X$-schéma étale muni d'une
action de $G$ compatible à celle définie sur $X$. Si le morphisme
$\varphi :U\rightarrow X$ est une immersion ouverte ou encore s'il
est affine, alors $F\mid U$ est adapté et $\widehat{H}^{*}_G(U;F\mid
U)=0$.
\end{cor}
\begin{proof}
On a le diagramme commutatif
\[ \xymatrix{
 U  \ar[d]_{s} \ar[r]^{\varphi} & X \ar[d]_{\pi}   \\
 U/G  \ar[r]_{\phi}& X/G
} \]
 Soient $t\in U$ et $x=\varphi(t)$. On a $I_t\subseteq I_x=\{1\}$, où $I_t$ et $I_x$
sont les groupes d'inerties respectivement aux points $t$ et $x$. Le
groupe $G$ opère sur $U$ sans inertie, donc $s$ est un revêtement
étale galoisien. Les morphismes $s$ et $\phi \circ s=\pi \circ
\varphi$ sont étales et $X/G$ est localement noethérien (car $X$
l'est). Il suit que $\phi$ est étale (cf \cite{SGA1} V.3.3). Pour
simplifier les notations, on pose $A:=\pi_G^*F$. On a $$F=\pi
^*A\,\, et\,\,\varphi^*F=(\pi \circ \varphi)^*A=(\phi \circ
s)^*A=s^*(\phi ^*A).$$ Comme $\phi$ est étale, on a bien $H^q(V;\phi
^*A)=0$ pour tout $q\geq n$ et pour tout $(U/G)$-schéma étale $V$.
Autrement dit, $\varphi^*F$ est adapté, car
$s_*^G(\varphi^*F)=\phi^*A$. D'autre part, $\varphi$ est de type
fini et $X$ est localement noethérien donc $U$ et $U/G$ le sont
aussi. On obtient $\widehat{H}^{*}_G(U;F\mid U)=0$ grâce à la
proposition \ref{net}.
\end{proof}

\begin{rem}\label{cite restreint}
Lorsque $Y=X/G$ est localement noethérien et quasi-séparé, les
résultats de cette section sont valables en remplaçant le cite étale
de $Y$ par le cite étale restreint de $Y$ (cf \cite{SGA4} VII.3.2).
En effet, la catégorie des faisceaux sur ce cite est équivalente à
celle des faisceaux sur le cite étale de $Y$ (qui est équivalente à
celle des $G$-faisceaux sur $X$). De plus, tous les schémas
envisagés dans cette section sont des objets du cite étale restreint
de $Y$, c'est à dire des $Y$-schémas étales de présentation finie
(cf \cite{SGA1} VIII.3.6).

En effet, dans la preuve de (\ref{net}) le morphisme $X\rightarrow
X/G=Y$ est fini (cf \cite{Raynaud} X corollaire 1). Il suit que
$X/S\rightarrow Y$ est de type fini et séparé (cf \cite{SGA1}
V.1.5). Ces deux morphismes sont donc de présentation finie
(quasi-séparés et de type fini). Dans la preuve de (\ref{NET}), le
morphisme étale $U\rightarrow X$ est soit affine soit une immersion
ouverte (et $X$ est localement noethérien), il est donc de
présentation finie. De même, par (cf \cite{SGA1} VIII.3.6), le
morphisme $U/G\rightarrow Y$ est de présentation finie.

Soient $X\rightarrow Y$ un revêtement étale galoisien de groupe $G$
et $F$ un $G$-faisceau sur $X$. Lorsque $Y$ est localement
noethérien et quasi-séparé, on dira parfois que $F$ est adapté sur
$X$ s'il est adapté en tant que faisceau sur le cite étale restreint
de $Y$. Alors $\widehat{H}^{*}_G(U;F\mid U)=0$ pour tout morphisme
équivariant $U\rightarrow X$ qui est soit étale et affine soit une
immersion ouverte.
\end{rem}

\subsection{Localisation à un voisinage étale des points ramifiés.}

On suppose que $X$ est connexe et localement noethérien. Soit de
plus un groupe fini $G$ opérant fidèlement sur $X$ et de manière
admissible. On considère le morphisme $\pi:X\rightarrow X/G$. Soit
$Z$ l'ensemble des points de $X$ dont les groupes d'inertie sont non
triviaux. D'après (\cite{Raynaud} X corollaire 1), $Z$ est aussi
l'ensemble des points en lesquels $\pi$ est ramifié. On en déduit
que $Z$ est fermé dans $X$ (cf \cite{SGA1} I.3.3). On note $X'$ le
complémentaire ouvert de $Z$ dans $X$. Le sous-schéma $X'$ est
localement noethérien, stable par $G$ et l'action de $G$ sur $X'$
est admissible.

 Un voisinage étale de $Z$ dans $X$ est un morphisme étale affine $\varphi:U\rightarrow X$ qui est un isomorphisme au-dessus de $Z$ (i.e $U\times_{X}Z\rightarrow Z$ est un isomorphisme).

\begin{defn}
On appelle \emph{$G$-voisinage étale de $Z$ dans $X$} un voisinage
étale $\varphi:U\rightarrow X$ muni d'une action de $G$ compatible à
celle définie sur $X$, de sorte que $\varphi^{-1}(Z)$ soit
d'intersection non vide avec chaque composante connexe de $U$.
\end{defn}

Si $U$ peut être muni d'une action de $G$ qui fait de
$\varphi:U\rightarrow X$ un $G$-voisinage étale de $Z$, alors cette
action est définie de manière unique. En effet, soit $\sigma\in G$
et supposons qu'il existe $f,g:U\rightarrow U$ rendant le diagramme
 \[ \xymatrix{
  U  \ar[d]_{\varphi} \ar[r]^{f ; g} & U \ar[d]_{\varphi}   \\
 X  \ar[r]^{\sigma}& X
} \]
  commutatif. Soient $p\in Z$ , $q$ et $q'$ les points de $U$ au-dessus de
$p$ et $\sigma(p)$ respectivement et $V$ la composante connexe de
$U$ contenant $q$. Alors $f(q)=g(q)=q'$ et $g^{-1}\circ f(q)=q$.
Comme $k(q)\simeq k(p)$,  $g^{-1}\circ f$ induit l'identité sur
$k(q)$. On en déduit que $g^{-1}\circ f\mid_{V}=Id_{V}$ (\cite
{Milne} I.3.13), puis que $ g^{-1}\circ f$ induit l'identité sur
chaque composante connexe de $U$.

\begin{thm} \label{localisation}
Soit $\varphi:U\rightarrow X$ un $G$-voisinage étale de $Z$ dans $X$
et $F$ un $G$-faisceau sur $X$. On suppose que $F\mid X'$ est
adapté.
 Alors il y a un isomorphisme canonique $$\widehat{H}_{G}^{*}(X;F)\simeq\widehat{H}_{G}^{*}(U;\varphi^{*}F).$$
\end{thm}

\begin{proof}
On remarque d'abord que $\varphi:U\rightarrow X$ se factorise à
travers $j:V\rightarrow X$, où $V:=Im(\varphi)$ est un ouvert de $X$
(pour la topologie de Zariski) et $j$ l'inclusion. On note
$\psi:U\rightarrow V$ le morphisme satisfaisant
$\varphi=j\circ\psi$. Alors $G$ opère sur $V$, $U$ et $X$ de manière
compatible. Autrement dit, le diagramme suivant est commutatif:

\[ \xymatrix{
  G\times U  \ar[d]_{1_G \times \psi} \ar[r] & U \ar[d]_{\psi}   \\
 G\times V \ar[d]_ {1_G \times j}\ar[r]& V \ar[d]_j \\
 G \times X \ar[r]                & X
} \]
 En particulier $j^{*}F$ et $\varphi^{*}F$ sont des $G$-faisceaux
sur $V$ et $U$ respectivement. D'autre part, le morphisme $\varphi$
est affine donc $\psi$ et $j$ le sont aussi. On note
$U':=X'\times_{X}U$, $V':=X'\times_{X}V$ et on commence par montrer
le résultat suivant:

{\bf \emph{(a) Les groupes $\widehat{H}_{G}^{*}(X;F)$ et
$\widehat{H}_{G}^{*}(V;j^{*}F)$ sont canoniquement isomorphes.}}

En posant $Z^{n}=Z^{n}(F)$, on a
$C^{n}(F)(V)=C^{0}(Z^{n})(V)=\Prod_{j\circ\beta\in
M}(Z^{n})_{\beta}$. La résolution de Godement (non tronquée) donne
la suite exacte de complexes de $\mathbb{Z}[G]$-modules
$$C^{*}(F)(X)\rightarrow C^{*}(F)(V)\rightarrow0.$$
Soit $K^{*}$ le complexe de $\mathbb{Z}[G]$-modules qui rend exacte
la suite
$$0\rightarrow K^{*}\rightarrow C^{*}(F)(X)\rightarrow C^{*}(F)(V)\rightarrow0.$$

Supposons avoir montré que  $H^{*}(Tot(W_*;K^*))=0$. Comme le
foncteur $j^{*}$ est exact et qu'il préserve les faisceaux flasques,
$$0\rightarrow j^{*}F\rightarrow
j^{*}C^{0}(F)\rightarrow...\rightarrow j^{*}C^{n}(F)\rightarrow...$$
est une résolution de $j^{*}F$ par des faisceaux flasques. Le fait
que $j:V\rightarrow X$ soit compatible à l'action de $G$ fait des
$j^{*}C^{n}(F)$ des $G$-faisceaux et des $j^{*}(d^{n})$ des
morphismes de $G$-faisceaux. On a donc une $G$-résolution flasque de
$j^{*}(F)$. D'autre part, $j^{*}(C^{n}(F))(V)=C^{n}(F)(V)$, donc la
cohomologie du complexe $C^{*}(F)(V)$ permet de calculer la
cohomologie équivariante du faisceau $j^{*}F$. Les $W_i$ sont des
$\mathbb{Z}[G]$-modules projectifs, donc la suite
$$0\rightarrow Tot(W_*;K^{*})\rightarrow Tot(W_*;C^{*}(F)(X))\rightarrow Tot(W_*;C^{*}(F)(V))\rightarrow0$$
est encore exacte. La suite exacte longue de cohomologie qui lui est
associée donne l'isomorphisme cherché.

Il suffit donc pour prouver \textit{(a)} de montrer que
$H^{*}(Tot(W_*;K^*))=0$. Soient $N\,:=\{\alpha\in M;\, Im(\alpha)\in
V\}$ et $\Delta:=\{\alpha\in M;\, Im(\alpha)\in Z\}$. On a alors
$$C^n(F)(X')=\Prod_{\alpha\in M-\Delta}(Z^{n})_{\alpha},\,\,
C^n(F)(V')=\Prod_{\alpha\in N-\Delta}(Z^{n})_{\alpha},$$
$$et\,\,\,K^{n}=\Prod_{\alpha\in M-N}(Z^{n})_{\alpha}=\Prod_{\alpha\in (M-\Delta)-(N-\Delta)}(Z^{n})_{\alpha}.$$
On en tire la suite exacte $$0\rightarrow K^{*}\rightarrow
C^{*}(F)(X')\rightarrow C^{*}(F)(V')\rightarrow0.$$ Les $W_i$ sont
projectifs et la suite
$$0\rightarrow Tot(W_{*};K^{*})\rightarrow Tot(W_{*};C^{*}(F)(X'))\rightarrow Tot(W_{*};C^{*}(F)(V'))\rightarrow0$$ est exacte. La suite exacte longue de cohomologie associée devient donc $$...\rightarrow H^q(Tot(W_*;K^*))\rightarrow \widehat{H}_{G}^{q}(X';F)\rightarrow \widehat{H}_{G}^{q}(V';F)\rightarrow...$$
$X'$ est localement nothérien et $G$ opère sur $X'$ de façon
admissible et sans inertie. De plus le morphisme $V'\rightarrow X'$
est une immersion ouverte compatible à l'action de $G$. Le
corollaire \ref{NET} s'applique, d'où $$\widehat{H}_{G}^{q}(X';F)=
\widehat{H}_{G}^{q}(V';F)=0,$$ et $H^{*}(Tot(W_*;K^*))=0$.

{\bf \emph{(b) Les groupes $\widehat{H}_{G}^{*}(V;F)$ et
$\widehat{H}_{G}^{*}(U;\psi^{*}F)$ sont canoniquement isomorphes.}}
Pour simplifier les notations, on désigne par $F$ le $G$-faisceau
$j^*F$ sur $V$. De plus, $M$ est maintenant un $G$-système de points
géométriques sur $V$, et $\Delta$ est toujours l'ensemble des points
géométriques de $M$ dont l'image est un point de $Z$.

Le morphisme $\psi:U\rightarrow V$ est surjectif, donc pour tout
$\alpha$ élément de $M$, il existe au moins un point géométrique de
$U$ au-dessus de $\alpha$. En effet, soient $v$ l'image de $\alpha$
et $U_v=U\times_{V}Spec(k(v))$ la fibre de $\psi$ au-dessus de $v$.
Alors $U_v=Spec(A)$ où $A$ est une $k(v)$-algèbre étale non nulle,
car $\psi$ est étale et surjectif. Si l'on note $r$ le degré de $A$
sur $k(v)$, il existe exactement $r$ points géométriques de $U$
au-dessus de $\alpha$.

 On a donc la suite exacte $$0\rightarrow C^{*}(F)(V)\rightarrow C^{*}(F)(U).$$
En posant $I^{n}:=C^{n}(F)(U)/C^{n}(F)(V)$, on obtient la suite
exacte
$$0\rightarrow C^{*}(F)(V)\rightarrow C^{*}(F)(U)\rightarrow I^{*}\rightarrow0.$$

On suppose là aussi avoir montré que $H^{*}(Tot(W_{*};I^{*}))=0$.
Comme ci-dessus,
$$0\rightarrow\psi^{*}F\rightarrow\psi^{*}C^{0}(F)\rightarrow...\rightarrow\psi^{*}C^{n}(F)\rightarrow...$$
est une $G$-résolution flasque de $\psi^{*}F$, et le complexe
$C^{*}(F)(U)$ permet de calculer la cohomologie équivariante du
faisceau $\psi^{*}F$. On a de nouveau la suite exacte
$$0\rightarrow Tot(W_*;C^{*}(F)(V))\rightarrow Tot(W_*;C^{*}(F)(U))\rightarrow Tot(W_*;I^{*})\rightarrow0,$$
 et la suite exacte longue de cohomologie associée donne l'isomorphisme voulu.

 Il reste à montrer que $H^{*}(Tot(W_{*};I^{*}))=0$.
 Si $\alpha\in \Delta$, $U$ possède un et un seul point géométrique au-dessus
de $\alpha$, car $U$ est un voisinage étale de $Z$.

D'où $C^{n}(F)(V)=(\Prod_{\alpha\in
\Delta}Z_{\alpha}^{n})\times(\Prod_{\alpha\in
M-\Delta}Z_{\alpha}^{n})$ , $C^{n}(F)(U)=(\Prod_{\alpha\in
\Delta}Z_{\alpha}^{n})\times(\Prod_{\psi\circ\beta\in M-\Delta}\,\,
Z_{\beta}^{n})$.
$$et\,\,I^{n}:=C^{n}(F)(U)/C^{n}(F)(V)=\prod_{\alpha\in M-\Delta}((\Prod_{\psi\circ\beta=\alpha}Z_{\beta}^{n})/Z_{\alpha}^{n})=C^{n}(F)(U')/C^{n}(F)(V'),$$
où $Z_{\alpha}^{n}$ se plonge diagonalement dans
$\Prod_{\psi\circ\beta=\alpha}Z_{\beta}^{n}$. On en déduit la suite
exacte
 $$0\rightarrow C^{*}(F)(V')\rightarrow C^{*}(F)(U')\rightarrow I^{*}\rightarrow0,$$
et enfin $$0\rightarrow Tot(W_*;C^{*}(F)(V'))\rightarrow
Tot(W_*;C^{*}(F)(U'))\rightarrow Tot(W_*;I^{*})\rightarrow0.$$ La
suite exacte longue associée à cette dernière devient
$$...\rightarrow \widehat{H}_{G}^{q}(V';F)\rightarrow \widehat{H}_{G}^{q}(U';F)\rightarrow H^q(Tot(W_*;I^{*}))\rightarrow...$$
Le schéma $X'$ est localement nothérien et $G$ opère sur $X'$ de
façon admissible et sans inertie. Le morphisme $V'\rightarrow X'$
est l'immersion d'un ouvert de $X'$ stable sous l'action $G$ et le
morphisme $U'\rightarrow X'$ est affine (car $\varphi$ l'est), étale
et compatible à l'action de $G$. Le corollaire \ref{NET} s'applique
et montre que
$\widehat{H}_{G}^{*}(V';F)=\widehat{H}_{G}^{*}(U';F)=0$. On obtient
bien $H^{*}(Tot(W_{*};I^{*}))=0$, ce qui achève la preuve du
théorème.

\end{proof}

\subsection{Hensélisation.}
On suppose que $X$ est un schéma connexe noethérien sur lequel le
groupe fini $G$ opère fidèlement. On conserve les hypothèses et les
notations de la section 2.2. Soit $X$ un espace topologique muni
d'une action d'un groupe fini $G$ et $Z$ un sous-espace fermé stable
par $G$. Si $U$ est un ouvert de $X$ contenant $Z$, alors $\bigcap
g(U)$ est un ouvert stable sous l'action de $G$ contenant $Z$ et
contenu dans $U$. On procède de la même manière pour montrer le
lemme suivant.

\begin{lem} L'ensemble des $G$-voisinages étales de $Z$ dans $X$ forme un système cofinal dans celui des voisinages étales de $Z$.
\end{lem}

\begin{proof}
Soit $\varphi:U\rightarrow X$ un voisinage étale de $Z$ dans $X$. On
note $_{\sigma}U$ le $X$-schéma étale $\sigma\circ
\varphi:U\rightarrow X$ et on pose $W:=\Prod_{\sigma\in
G}(_{\sigma}U)$, le produit fibré étant pris sur $X$. Le groupe $G$
opère sur lui-même par translations. Ceci définit une action de $G$
sur $W$ par permutation des coordonnées qui est compatible à
l'action de $G$ sur $X$. Chacun des $_{\sigma}U$ est un voisinage
étale de $Z$ dans $X$. Il suit que $W\times_{X}Z\rightarrow Z$ est
un isomorphisme. D'autre part, $W\rightarrow X$ est affine puisque
tous les $_{\sigma}U\rightarrow X$ le sont. On note $GU$ le
sous-schéma de $W$ formé des composantes connexes de $W$ contenant
au moins un point s'envoyant dans $Z$. Alors $GU$ est un
$G$-voisinage étale de $Z$ dans $X$ de sorte qu' il y a un morphisme
canonique de $X$-schémas $GU\rightarrow U$. De plus, si
$f:V\rightarrow U$ est un morphisme de $X$-schémas, $f$ induit
canoniquement un morphisme $GV\rightarrow GU$ au-dessus de $X$ et
compatible à l'action de $G$.
\end{proof}

On note $\widetilde{Z}$ la limite projective des voisinages étales
de $Z$ dans $X$. D'après ce qui précède, c'est aussi la limite
projective des $G$-voisinages étales de $Z$ dans $X$. L'action de
$G$ sur ces derniers passe à la limite pour donner une action de $G$
sur $\widetilde{Z}$ compatible à celle définie sur $X$. D'autre
part, lorsque $X$ est noethérien, tout voisinage étale $U$ de $Z$
dans $X$ est noethérien (en particulier quasi-compact) puisque le
morphisme $U\rightarrow X$ est de type fini.

On note $i:Z\rightarrow X$, $s:Z\rightarrow \widetilde{Z}$ et
$\phi:\widetilde{Z}\rightarrow X$ les morphismes canoniques. Ces
morphismes commutent à l'action de $G$ et satisfont $i=\phi \circ
s$.


\begin{thm} \label{inductive}
Soit $F$ un $G$-faisceau sur $X$ tel que $F\mid X'$ soit adapté.
Alors le morphisme canonique $\widehat{H}_{G}^{*}(X;F) \rightarrow
\widehat{H}_{G}^{*}(\widetilde{Z};\phi^*F)$ est un isomorphisme.
\end{thm}

\begin{proof} On note $0\rightarrow F\rightarrow C^0\rightarrow ...$ la résolution de Godement de $F$. On a les égalités suivantes (où toutes les limites inductives sont prises sur l'ensemble des $G$-voisinages étales de $Z$ dans $X$):

\begin{eqnarray}
\label{un}\widehat{H}_{G}^{*}(X;F)&=& \displaystyle{\lim_{\longrightarrow }}\,\,  \widehat{H}_{G}^{*}(U;F)\\
\label{deux}&=&H^*(\displaystyle{\lim_{\longrightarrow }}\,\, Tot^*(W_*;C^*(U)))\\
\label{trois}&=&H^*\Big(\bigoplus_{i+j=*}\displaystyle{\lim_{\longrightarrow}}\,\,Hom_{\mathbb{Z}[G]}(W_i;C^j(U))\Big)\\
\label{quatre}&=&H^*\Big(\bigoplus_{i+j=*} Hom_{\mathbb{Z}[G]}(W_i;\displaystyle{\lim_{\longrightarrow }}\,\,C^j(U))\Big)\\
\label{cinq}&=&H^*\Big(\bigoplus_{i+j=*} Hom_{\mathbb{Z}[G]}(W_i;\phi^*C^j(F)(\widetilde{Z}))\Big)\\
\label{six}&=&\widehat{H}_{G}^{*}(\widetilde{Z};\phi^*F)
\end{eqnarray}

L'égalité \eqref{un} est vraie car tous les morphismes
$\widehat{H}_{G}^{*}(U;F)\rightarrow\widehat{H}_{G}^{*}(V;F)$ sont
des isomorphismes. Etant donné un système inductif de complexes de
${\mathbb{Z}[G]}$-modules, la limite inductive commute avec la
cohomologie, ce qui donne \eqref{deux}. L'égalité \eqref{trois} est
vraie car les limites inductives commutent avec les sommes directes.
L'égalité \eqref{quatre} est vraie car les $W_i$ sont des
${\mathbb{Z}[G]}$-modules libres de type fini. Le foncteur $\phi^*$
est exact, donc $0\rightarrow \phi^*F\rightarrow
\phi^*C^0\rightarrow ...$ est une résolution de $\phi^*F$. Les
voisinages étales sont quasi-compacts et les morphismes de
transition sont affines, donc les morphismes
$\displaystyle{\lim_{\longrightarrow }}\,\,H^q(U;C^j)\rightarrow
H^q(\widetilde{Z};\phi^*C^j)$ sont bijectifs. Ainsi, les faisceaux
$\phi^*C^j$ sont acycliques pour le foncteur des sections globales .
Les égalités \eqref{cinq} et \eqref{six} s'en déduisent.

\end{proof}

Supposons que $Z$ soit contenu dans un ouvert affine $Spec(A)$. Le
sous-schéma fermé $Z$ est alors défini par un idéal $I=\sqrt{I}$ de
$A$ de sorte que $Z=V(I)$ et $Z\simeq Spec(A/I)$. Si
$(\widetilde{A};\widetilde{I})$ désigne l'hensélisé du couple
$(A;I)$ (cf \cite {Raynaud} XI), alors
$\widetilde{Z}:=Spec(\widetilde{A})$.

\begin{cor} \label{analogue}
On conserve les hypothèses du théorème \ref{inductive}. Lorsque $F$
est de torsion et que $Z$ est contenu dans un ouvert affine, le
morphisme canonique $$\widehat{H}_{G}^{*}(X;F) \rightarrow
\widehat{H}_{G}^{*}(Z;i^*F)$$ est un isomorphisme.
\end{cor}

\begin{proof} Soit $s:Z\rightarrow \widetilde{Z}$ l'immersion fermé
canonique. Comme $Z$ est contenu dans un ouvert affine,
$\widetilde{Z}$ est affine. De plus, le $G$-faisceau $\phi^*F$ sur
$\widetilde{Z}$ est de torsion lorsque $F$ l'est. Soient $I^*$ et
$J^*$ des $G$-résolutions injectives de $\phi^*F$ et
$s^*\phi^*F=i^*F$ respectivement. Le morphisme canonique de
complexes de $\mathbb{Z}[G]$-modules $I^*(\widetilde{Z})\rightarrow
J^*(Z)$ induit un morphisme de suites spectrales et des
isomorphismes sur les groupes de cohomologie étale (cf \cite{Huber}
0.1). Comme les suites spectrales en question convergent
respectivement vers $\widehat{H}_{G}^{*}(\widetilde{Z};F)$ et
$\widehat{H}_{G}^{*}(Z;s^*F)$, ces deux groupes (gradués) sont
isomorphes. En utilisant le théorème \ref{inductive}, on obtient le
résultat annoncé.
\end{proof}

\begin{rem}\label{Z} Soit $F$ un $G$-faisceau quelconque sur $X$ qui est adapté sur
$X'$. Pour appliquer la preuve précédente à $F$, il suffit de
montrer que les morphismes canoniques
$H^q(\widetilde{Z};\phi^*F)\rightarrow H^q(Z;i^*F)$ sont bijectifs.
Cette condition est toujours vérifiée lorsque $\widetilde{Z}$ est
une somme finie de spectres d'anneaux locaux henséliens (cf
\cite{SGA4} VIII 8.6.).
\end{rem}


Soit $F$ un faisceau de torsion sur un schéma $X$. Pour tout nombre
premier $l$, on note $F(l)$ le sous-faisceau de $l$-torsion de $F$.
Le morphisme canonique $\bigoplus_lF(l)\rightarrow F$ est un
isomorphisme.

\begin{prop}
Soit $G$ un groupe fini opérant fidèlement et de manière admissible
sur $X$. Si $X/G$ est un ouvert du spectre de l'anneau d'entiers
d'un corps de nombres totalement imaginaire, alors tout $G$-faisceau
sur $X$ est adapté sur $X'$.

Soit $X$ de type fini et séparé sur $Spec(\mathbb{Z})$, et soit $F$
un $G$-faisceau de torsion sur $X$. Si $F(2)={0}$ ou si aucun corps
résiduel de $X$ n'est ordonnable, alors $F$ est adapté sur $X'$.

Soit $X$ un schéma de type fini séparé sur un corps $k$, sur lequel
$G$ opère par $k$-automorphismes. Si $k$ est de $l$-dimension
cohomologique finie, tout $G$-faisceau de $l$-torsion sur $X$ est
adapté sur $X'$. Si $k$ est de dimension cohomologique finie, tout
$G$-faisceau de torsion sur $X$ est adapté sur $X'$.
\end{prop}

\begin{proof}
Si $X/G$ est un ouvert du spectre de l'anneau d'entiers d'un corps
de nombres totalement imaginaire, tout $X'/G$-schéma étale (connexe)
est aussi un ouvert du spectre de l'anneau d'entiers d'un corps de
nombres totalement imaginaire. Dans ce cas, le résultat est une
conséquence de \cite{Deninger} 4.6.

On note $S$ la base $Spec(\mathbb{Z})$ ou $Spec(k)$. Soient $X$ de
type fini et séparé sur $S$, et $F$ comme dans l'énoncé. L'immersion
ouverte $X'\rightarrow X$ est séparée et de type fini, car $X$ est
noethérien. Ainsi $X'$ est séparé et de type fini sur $S$. D'après
(\cite{SGA1} V.1.5), le schéma $X'/G$ est aussi de type fini et
séparé sur $S$.

On pose $\pi:X'\rightarrow X'/G$ et $A:=\pi_*^G F$. Le faisceau $A$
est de torsion (respectivement de $l$-torsion) si $F$ l'est. D'après
la remarque (\ref{cite restreint}), il suffit de montrer qu'il
existe un entier $n$ de sorte que $H^q(U;A)=0$ pour tout $q\geq n+1$
et pour tout $U$ étale de présentation finie sur $X'/G$. Dans ce
cas, le schéma $U$ est quasi-compact et quasi-séparé, donc la
cohomologie étale de $U$ commute aux limites inductives filtrantes
de faisceaux (cf \cite{SGA4} VII.3.2).

On pose $n:=2\,dim(X)+1$ pour $S=Spec(\mathbb{Z})$ et $n:=2\,dim(X)+
cd_l(k)$ (respectivement $n:=2\,dim(X)+ cd(k)$) pour $S=Spec(k)$. On
obtient, pour tout $q\geq n+1$,
$$H^q(U;A)=\bigoplus_l H^q(U;A(l))=0.$$
En effet, le $S$-schéma $U$ est de type fini, donc la deuxième
égalité est vraie grâce à (\cite{SGA4} X.6.2) et (\cite{Milne}
VI.1.4). Ainsi, le $G$-faisceau $F$ est adapté sur $X'$.

\end{proof}

\section{Application à la topologie arithmétique.}

Soit $X$ le spectre de l'anneau d'entiers $D$ d'un corps de nombres
$L$ sur lequel un groupe fini $G$ opère. On pose $K:=L^{G}$,
$Y=X/G=Spec(D^G)$ et on considère l'extension galoisienne $L/K$. On
suppose $L$ et $K$ totalement imaginaires. Les groupes de
cohomologie étale du faisceau du groupe multiplicatif sur le site
$X_{et}$ sont donnés ci-dessous (cf \cite{Ma}).
$$
\begin{array}{rcl}
H^{q}(X;\mathbb{G}_m)&=& U_L\,\,\,\, pour\,\, q=0,\\
          &=& Cl(L)\,\,\,\, pour\,\, q=1,\\
&=&0\,\,\,\, pour\,\, q=2,\\
&=&\mathbb{Q}/\mathbb{Z}\,\,\,\, pour\,\, q=3,\\
&=&0\,\,\,\, pour\,\, q\geq4,\\
\end{array}$$
où $U_L$ désigne le groupe des unités de $L$, $\mu$ le groupe
cyclique des racines de l'unité et $Cl(L)$ le groupe des classes.
L'opération de $G$ sur les groupes $U_L$ et $Cl(L)$ (en tant que
groupes de cohomologie du $G$-faisceau $\mathbb{G}_m$) est l'action
naturelle. De plus, le groupe de Galois opère trivialement sur
$H^{3}(X;\mathbb{G}_m)=\mathbb{Q}/\mathbb{Z}$.

 Le morphisme canonique $\pi :X\rightarrow Y$ est fini. On note $Z$ le sous-schéma fermé de $X$ constitué des points en lesquels $\pi$ est ramifié, $X'$ son complémentaire ouvert et $s$ le cardinal de $Z$.

\subsection{Calcul des groupes $\widehat{H}_{G}^{*}(X;\mathbb{G}_{m})$.}
Le quotient $X'/G$ est un ouvert de $Y$. Soit $V$ un $X'/G$-schéma
étale. Alors $V$ est un ouvert du spectre de l'anneau d'entiers d'un
corps de nombres totalement imaginaire et on a
$H^q(V;\mathbb{G}_m)=0$ pour tout $q \geq 4$ (cf \cite{Milne2}
2.2.1). Le faisceau $\mathbb{G}_m\mid X'$ est adapté donc le
théorème de localisation s'applique dans ce cas.

  Soient $\mathfrak{a}$ le produit des idéaux $(\mathfrak{p}_{i})_{1\leq i\leq s}$ de $D$ ramifiés
dans l'extension $L/K$ et $Z=V(\mathfrak{a})$. Soient
$(\widetilde{D};\widetilde{\mathfrak{a}})$ l'hensélisé du couple
$(D;\mathfrak{a})$ et $\widetilde{Z}:=Spec(\widetilde{D})$ (cf
\cite{Raynaud}). Alors $\widetilde{D}=\Prod_{1\leq i\leq
s}\widetilde{D}_{\mathfrak{p}_{i}}$, où
$\widetilde{D}_{\mathfrak{p}_{i}}$ est l'anneau local constitué des
éléments de la complétion de $D$ pour la valeur absolue
$\mathfrak{p}_{i}$-adique qui sont algébriques sur $\mathbb{Q}$.
C'est un anneau de valuation discrète hensélien à corps résiduel
fini. On note $U_{\mathfrak{p}_i}$ le groupe des unités de
$\widetilde{D}_{\mathfrak{p}_{i}}$ et $\widetilde{Z}_{i}$ le spectre
de $\widetilde{D}_{\mathfrak{p}_{i}}$.

Le théorème de localisation donne l'isomorphisme
$\widehat{H}_{G}^{*}(X;\mathbb{G}_{m})\simeq
\widehat{H}_{G}^{*}(\widetilde{Z};\mathbb{G}_{m})$. D'autre part, on
a $\widetilde{Z}=\Coprod \widetilde{Z}_{i}$ et
$H^{q}(\widetilde{Z};\mathbb{G}_{m})=\Prod
H^{q}(\widetilde{Z}_{i};\mathbb{G}_{m})=\Prod U_{\mathfrak{p}_i}$
pour $q=0$ et $0$ sinon. En effet,
$H^q(\widetilde{Z}_i;\mathbb{G}_m)\simeq
H^q(Spec(D/\mathfrak{p}_i);\mathbb{G}_m)$ pour tout $q\geq 1$ (cf
\cite{Strano} 4.1), or ces groupes sont nuls car un corps fini est
C1 et donc de dimension cohomologique $\leq 1$ (cf \cite {Serre2}).
La suite spectrale de cohomologie équivariante de $\mathbb{G}_m$ sur
$\widetilde{Z}$ dégénère et on obtient
\begin{equation}\label{huit}
\widehat{H}_{G}^{*}(X;\mathbb{G}_{m})\simeq
\widehat{H}_{G}^{*}(\widetilde{Z};\mathbb{G}_{m})=\widehat{H}^{*}(G;\Prod_{1\leq
i\leq s}U_{\mathfrak{p}_i}).
\end{equation}

Si $\mathfrak{q}$ est un idéal premier de $K$, on pose
$U^{\mathfrak{q}}:=\Prod_{\mathfrak{p}\mid
\mathfrak{q}}U_{\mathfrak{p}}$ et on a la décomposition en
$\mathbb{Z}[G]$-modules homogènes $\Prod_{1\leq i\leq
s}U_{\mathfrak{p}_i}=\Prod_{\mathfrak{q}\in\Omega}U^{\mathfrak{q}}$,
où $\Omega$ désigne l'ensemble des idéaux premiers de $K$ qui se
ramifient dans $L$. Lorsque $\mathfrak{p}$ est un idéal premier de
$L$, on note $G_{\mathfrak{p}}$ le groupe de décomposition en
l'idéal $\mathfrak{p}$. Alors $U^{\mathfrak{q}}$ est le
$\mathbb{Z}[G]$-module induit
$M_{G}^{G_{\mathfrak{p}}}(U_{\mathfrak{p}})$, où l'on a choisi un
idéal $\mathfrak{p}$ premier de $L$ divisant $\mathfrak{q}$.
Finalement, on a l'identification
\begin{equation}
\label{neuf} \widehat{H}^{p}(G;\Prod_{1\leq i\leq
s}U_{\mathfrak{p}_i})=\Prod_{\mathfrak{q}\in\Omega}{\widehat{H}^{p}(G_{\mathfrak{p}};U_{\mathfrak{p}})},
\end{equation}
où l'on a choisi, pour tout $\mathfrak{q}$ dans $\Omega$, un idéal
$\mathfrak{p}$ premier de $L$ divisant $\mathfrak{q}$.

Les relations (\ref {huit}) et (\ref {neuf}) montrent le théorème
suivant.
\begin{thm}
On a l'isomorphisme $$\widehat{H}^{n}_G(X;\mathbb{G}_m)\simeq\Prod
\widehat{H}^{n}(G_{\mathfrak{p}};U_{\mathfrak{p}}),$$ où le produit
est pris sur tous les idéaux premiers non nuls de $K$.
\end{thm}

La théorie du corps de classe local (cf \cite {Milne2} I Appendix A)
montre que le groupe
$\widehat{H}^{0}(G_{\mathfrak{p}};U_{\mathfrak{p}})$ est isomorphe
au sous-groupe d'inertie $I_{\mathfrak{q}}^{(ab)}$ de l'abélianisé
du groupe de décomposition $G_{\mathfrak{p}}$, c'est à dire le noyau
du morphisme $G_{\mathfrak{p}}^{ab}\rightarrow
Gal(k(\mathfrak{p})/k(\mathfrak{q}))$.

Soit $B$ l'anneau d'entiers de $K$ et $\widetilde{B}_{\mathfrak{q}}$
l'anneau de valuation discrète hensélien obtenu comme ci-dessus. On
note $U_{\mathfrak{q}}$ le groupe des unités de
$\widetilde{B}_{\mathfrak{q}}$, $L_{\mathfrak{p}}$ et
$K_{\mathfrak{q}}$ les corps de fractions de
$\widetilde{D}_{\mathfrak{p}}$ et $\widetilde{B}_{\mathfrak{q}}$. On
considère la suite exacte de $\mathbb{Z}[G_{\mathfrak{p}}]$-modules
$$0\rightarrow U_{\mathfrak{p}}\rightarrow
L_{\mathfrak{p}}^{\times}\rightarrow\mathbb{Z}\rightarrow0$$ où la
flèche $L_{\mathfrak{p}}^{\times}\rightarrow\mathbb{Z}$ est donnée
par la valuation $\mathfrak{p}$-adique. En utilisant le fait que
l'extension $L_{\mathfrak{p}}/K_{\mathfrak{q}}$ est galoisienne et
le théorème de Hilbert 90, la suite exacte longue de cohomologie
donne la suite exacte $$0\rightarrow U_{\mathfrak{q}}\rightarrow
K_{\mathfrak{q}}^{\times}\rightarrow\mathbb{Z}\rightarrow\widehat{H}^{1}(G_{\mathfrak{p}};U_{\mathfrak{p}})\rightarrow0.$$
La restriction de la valuation $\mathfrak{p}$-adique à
$K_{\mathfrak{q}}^{\times}$ a pour image
$e_{\mathfrak{q}}\mathbb{Z}$, où $e_{\mathfrak{q}}$ désigne l'indice
de ramification associé à $\mathfrak{q}$ ($L/K$ est galoisienne,
donc cet entier ne dépend que de $\mathfrak{q}$). La dernière suite
exacte permet donc d'identifier les groupes
$\widehat{H}^{1}(G_{\mathfrak{p}};U_{\mathfrak{p}})$ et
$\mathbb{Z}/e_{\mathfrak{q}}\mathbb{Z}$. Le corollaire suivant
résume ces résultats.

\begin{cor} \label{calcul} On a les identifications suivantes, où les trois produits sont pris sur tous les idéaux premiers (non nuls) de $K$:
\begin{itemize}
\item $\widehat{H}^{0}_G(X;\mathbb{G}_m)= \Prod I_{\mathfrak{q}}^{(ab)}$.
\item $\widehat{H}^{1}_G(X;\mathbb{G}_m)= \Prod \mathbb{Z}/e_{\mathfrak{q}}\mathbb{Z}.$
\item $\widehat{H}^{n}_G(X;\mathbb{G}_m)= \Prod \mathbb{Z}/e_{\mathfrak{q}}\mathbb{Z}$ pour tout $n$ si $G$ est cyclique.
\end{itemize}
\end{cor}

\subsection{Etude de la suite spectrale de cohomologie équivariante.}

On conserve les mêmes notations. La deuxième page de la suite
spectrale relative au $G$-faisceau $\mathbb{G}_m$ est la suivante :

{\small{
$$
\begin{array}{cccccc}

\widehat{H}^{-3}(G;\mathbb{Q}/\mathbb{Z})&\widehat{H}^{-2}(G;\mathbb{Q}/\mathbb{Z})      &\widehat{H}^{-1}(G;\mathbb{Q}/\mathbb{Z})                     &\widehat{H}^{0}(G;\mathbb{Q}/\mathbb{Z})          &\widehat{H}^{1}(G;\mathbb{Q}/\mathbb{Z})&\widehat{H}^{-2}(G;\mathbb{Q}/\mathbb{Z})  \\[0.5cm]
0&0                     &0                     &0                     &0&0             \\[0.5cm]
\widehat{H}^{-3}(G;Cl(L))&\widehat{H}^{-2}(G;Cl(L))&\widehat{H}^{-1}(G;Cl(L))&\widehat{H}^0(G;Cl(L)) &\widehat{H}^1(G;Cl(L))    &\widehat{H}^2(G;Cl(L))  \\[0.5cm]
\widehat{H}^{-3}(G;U_L)&\widehat{H}^{-2}(G;U_L)&\widehat{H}^{-1}(G;U_L)&\widehat{H}^{0}(G;U_L) &\widehat{H}^1(G;U_L) &\widehat{H}^2(G;U_L)  \\[0.5cm]
\end{array}
$$
}} On a l'indentification
$\widehat{H}^{n}(G;\mathbb{Q}/\mathbb{Z})\simeq\widehat{H}^{n+1}(G;\mathbb{Z})$.
Lorsque $G$ est cyclique, la troisième ligne prend les valeurs $0$
pour les colonnes d'indice pair et
$\widehat{H}^{-2}(G;\mathbb{Z})\simeq G^{ab}=G$ pour celles d'indice
impair. On cherche des hypothèses sous lesquelles certaines
différentielles sont nulles.

On suppose désormais que le groupe $Cl(K)$ est trivial. Lorsque $G$
est cyclique, il existe une infinité d'idéaux premiers non nuls
inertes de $D$ (c'est une conséquence du théorème de densité de
Chebotarev (cf \cite{Sikora} théorème 2.3(3))). Cette condition est
d'ailleurs nécessaire. On peut alors choisir un point fermé $x$ de
$X$ (correspondant à un tel idéal $\mathfrak{p}$) fixé par $G$ et de
sorte que $\pi :X\rightarrow Y$ soit étale en $x$. Dans ces
conditions, $\mathfrak{p}$ est principal et on pose
$U=X-\{x\}=Spec(D_f)$, avec $\mathfrak{p}=D_f$. Les groupes de
cohomologie étale de $U$ à coefficients dans $\mathbb{G}_m$ sont les
suivants (cf \cite{Milne2} II.2.1):
$$
\begin{array}{rcl}
H^{q}(U;\mathbb{G}_m)&=& D_f^{\times}\,\,\,\, pour\,\, q=0,\\
          &=& Pic(U)=Cl(D_f)\,\,\,\, pour\,\, q=1,\\
&=&0\,\,\,\, pour\,\, q\geq2.\\
\end{array}$$

On désigne toujours par $\widetilde{Z}$ la limite projective des
voisinages étales de $Z$ dans $X$. Les $d_*^{**}$ sont les
différentielles de la suite spectrale  $E^{**}_*(X)$ de cohomologie
équivariante de $\mathbb{G}_m$ sur $X$. On note aussi
$E^{**}_*(\widetilde{Z})$ (respectivement $E^{**}_*(U)$) la suite
spectrale associée à $\mathbb{G}_m$ sur $\widetilde{Z}$
(respectivement sur $U$).
\begin{lem} \label{lem1}
Si le groupe $Cl(K)$ est trivial, la différentielle $d_2^{-1;1}$ est
nulle. Si de plus le groupe $G$ est cyclique, les différentielles
$d_2^{n;1}$ sont nulles lorsque $n$ est impair.
\end{lem}
\begin{proof}
On considère le morphisme de suites spectrales
$$h:E^{**}_*(X)\longrightarrow E^{**}_*(\widetilde{Z}).$$
On a en particulier le diagramme commutatif
\[ \xymatrix{
E_2^{-1;1}(X) \ar[d]_{d_2^{-1;1}} \ar[r]^{h_2^{-1;1}}  &E_2^{-1;1}(\widetilde{Z})=0 \ar[d]_{\delta_2^{-1;1}}   \\
E_2^{1;0}(X)\ar[r]^{h_2^{1;0}}         &E_2^{1;0}(\widetilde{Z}) }
\] Il suit que $h_2^{1;0}\circ d_2^{-1;1}=0$. D'autre part, la
flèche $h_2^{1;0}$ est donnée par le morphisme canonique
$\widehat{H}^{1}(G;U_L)\rightarrow \Prod_{\mathfrak{p}\in
Z}\widehat{H}^{1}(G;U_{\mathfrak{p}})$, qui est injectif car
$Cl(K)=0$ (cf lemme \ref {clef}). La deuxième affirmation suit en
utilisant la périodicité de la cohomolgie des groupes cycliques.
\end{proof}

\begin{lem}\label{lem2}
Si le groupe $Cl(K)$ est trivial et si $G$ est cyclique, toutes les
différentielles $d_3^{n;3}$ sont nulles.
\end{lem}
\begin{proof} Sous ces conditions, on dispose du morphisme de suites spectrales
$$j:E^{**}_*(X)\longrightarrow E^{**}_*(U).$$
Comme $\mathfrak{q}$ est principal, le morphisme $Cl(D)\rightarrow
Cl(D_f)$ est bijectif et induit un isomorphisme
$j_2^{n+3;1}:\widehat{H}^{n+3}(C_p;Cl(D))\rightarrow
\widehat{H}^{n+3}(C_p;Cl(D_f))$. Ce dernier induit à son tour (par
restriction aux noyaux des différentielles $d_2^{n+3;1}$ et
$\delta_2^{n+3;1}$) le morphisme $j_3^{n+3;1}$, qui est donc
injectif. Par ailleurs, le même argument que dans la preuve
précédente montre que $j_3^{n+3;1}\circ d_3^{n;3}$ est nul, ce qui
permet de conclure.
\end{proof}

\begin{lem} \label{lem3}
Si le groupe $Cl(K)$ est trivial, la différentielle $d_4^{-3;3}$ est
nulle. Si de plus $G$ est cyclique, toutes les différentielles
$d_4^{n;3}$ sont nulles.
\end{lem}
\begin{proof}
Il suffit de montrer que le morphisme $h_4^{1;0}$ est injectif.

On vérifie grâce au lemme \ref{lem1} que les groupes d'arrivée et de
départ du morphisme $h_2^{1;0}:E_2^{1;0}(X)\rightarrow
E_2^{1;0}(\widetilde{Z})$ restent inchangés jusqu'à la quatrième
page. Le morphisme $h_4^{1;0}$ s'identifie alors à $h_2^{1;0}$, qui
est injectif (cf lemme \ref {clef}), ce qui montre la première
affirmation.

On montre la deuxième en utilisant à nouveau la périodicité de la
cohomologie des groupes cycliques et l'égalité
$\widehat{H}^n(G;\mathbb{Q}/\mathbb{Z})=0$ lorsque n est pair.
\end{proof}

\subsection{Applications.}
Dans tout ce qui suit, $L/K$ est une extension galoisienne de corps
de nombres totalement imaginaires de groupe $G$. On note $X$
(respectivement $Y$) le spectre de l'anneau d'entiers de $L$
(respectivement de $K$).

\subsubsection{Majoration et minoration du nombre $s$ d'idéaux premiers ramifiés.}

On suppose que l'extension $L/K$ est cyclique d'ordre premier~$p$.
On note $C_p$ son groupe de Galois et $s$ le nombre d'idéaux
premiers (non nuls) ramifiés dans cette extension.

\begin{thm} \label{majoration} On la majoration
$ s \leq 1+dim_{\mathbb{F}_p}\widehat{H}^1(C_p;U_L/\mu)+
dim_{\mathbb{F}_p}\widehat{H}^0(C_p;Cl(L)).$
\end{thm}
\begin{proof}
Le théorème de localisation fournit les isomorphismes (cf corollaire
\ref{calcul})
$$\widehat{H}_{G}^{n}(X;\mathbb{G}_{m})\simeq
\widehat{H}_{G}^{n}(\widetilde{Z};\mathbb{G}_{m})\simeq\mathbb{F}_p^s.$$

 La suite spectrale converge vers $\mathbb{F}_p^s$ en tout degré, donc les groupes $\Prod_{i+j=1}E_{\infty}^{i;j}(X)$ et
$\mathbb{F}_p^s$ ont le même cardinal. Ainsi, l'inégalité
(immédiate) $\sharp (\Prod_{i+j=1}E_{\infty}^{i;j}(X)) \leq
\sharp(\Prod_{i+j=1}E_{2}^{i;j}(X))$ donne la majoration

\begin{equation}
\label{unmajoration}
s\leq dim_{\mathbb{F}_p}\widehat{H}^1(C_p;U_L)+
dim_{\mathbb{F}_p}\widehat{H}^0(C_p;Cl(L)).
\end{equation}
 Le groupe $\mu$ est cyclique, donc
$dim_{\mathbb{F}_p}\widehat{H}^1(C_p;\mu)\leq 1$. La suite exacte de
$C_p$-modules
$$0\rightarrow \mu \rightarrow U_L\rightarrow U_L/\mu\rightarrow0$$
donne la suite exacte $\widehat{H}^1(C_p;\mu)\rightarrow
\widehat{H}^1(C_p;U_L)\rightarrow \widehat{H}^1(C_p;U_L/\mu)$ et
l'inégalité
\begin{equation} \label{deuxmajoration}
dim_{\mathbb{F}_p}\widehat{H}^1(C_p;U_L)\leq
1+dim_{\mathbb{F}_p}\widehat{H}^1(C_p;U_L/\mu).
\end{equation}
Les relations (\ref{unmajoration}) et (\ref{deuxmajoration})
permettent de conclure.
\end{proof}

L'étude de la suite spectrale faite dans la section 4.2 est
nécessaire pour minorer le nombre~$s$.
\begin{thm} \label{minoration} Si le groupe $Cl(K)$ est trivial, alors $s\geq 1+dim_{\mathbb{F}_p}\widehat{H}^0(C_p;Cl(L)).$
\end{thm}
\begin{proof}
Les lemmes \ref{lem1}, \ref{lem2} et \ref{lem3} montrent les
égalités suivantes :
$$E_{\infty}^{-3;3}(X)=E_2^{-3;3}(X)=\mathbb{F}_p\,\, et \,\,E_{\infty}^{-1;1}(X)=E_2^{-1;1}(X)=\widehat{H}^{-1}(C_p;Cl(L)).$$
Cette suite spectrale converge vers
$\widehat{H}^n_{C_p}(X;\mathbb{G}_m)\simeq\mathbb{F}_p^s$ donc les
groupes $\Prod_{i+j=0}E_{\infty}^{i;j}(X)$ et $\mathbb{F}_p^s$ ont
donc le même cardinal. On obtient l'inégalité
$$ s\geq 1+dim_{\mathbb{F}_p}\widehat{H}^1(C_p;Cl(L)).$$
De plus, $\widehat{H}^1(C_p;Cl(L))$ et $\widehat{H}^0(C_p;Cl(L))$
ont le même cardinal car $Cl(L)$ est fini (\cite{Serre} VIII
proposition 8).
\end{proof}

\subsubsection{Lorsque $L$ possède une place finie stable par $G$.}
Soit alors $U=Spec(D_f)$ le complémentaire ouvert dans $X$ d'un
point fermé $x$ stable sous l'action de $G$. La suite spectrale
$\widehat{H}^p(G;H^q(U;\mathbb{G}_m))\Longrightarrow\widehat{H}^{p+q}_G(U;\mathbb{G}_m)$
ne possède que deux lignes non nulles. En observant son terme
initial, on obtient le résultat suivant.
\begin{prop}\label{complémentaire}
Si le groupe de Galois (quelquonque) de l'extension $L/K$ fixe une
place finie de $L$, on a la suite exacte longue

 $
...\xrightarrow{d_2^{n-2;1}} \widehat{H}^n(G;D_f^{\times})\rightarrow \widehat{H}_G^n(U;\mathbb{G}_m)\rightarrow \widehat{H}^{n-1}(G;Cl(D_f)) \\[2mm]
\mbox{}\hfill\xrightarrow{d_2^{n-1;1}} \widehat{H}^{n+1}(G;D_f^{\times})\rightarrow \widehat{H}_G^{n+1}(U;\mathbb{G}_m)\rightarrow \widehat{H}^{n}(G;Cl(D_f))\xrightarrow{d_2^{n;1}}  \\
$
\end{prop}
Si $L/K$ n'est ramifiée en aucun point de $U$, toutes les
différentielles $d_2^{n;1}$ sont des isomorphismes. Si maintenant
$G$ est cyclique, cette suite exacte s'exprime avec le produit des
groupes d'inertie pris sur l'ensemble des points fermés de $U/G$. De
plus, si $Cl(D_f^G)={0}$, le lemme \ref{lem1} s'applique et les
différentielles $d_2^{n;2}$ sont nulles pour $n$ impair. La suite
exacte longue se réduit à une suite exacte courte à six termes.

\subsubsection{Revêtements cycliques d'une sphère homologique.}
Dans le contexte de la topologie arithmétique, Niranjan Ramachandran
(cf \cite{Ramachandran}) a proposé la définition suivante.

\begin{defn}
Le spectre $Y$ de l'anneau d'entiers d'un corps de nombre $K$ est
une \emph{3-sphère à homologie entière} si $H^p(Y;\mathbb{G}_m)=0$
pour $p\neq0,3$ et si $H^0(Y;\mathbb{G}_m)$ est de torsion. $Y$ est
une 3-sphère à homologie rationnelle si $H^0(Y;\mathbb{G}_m)$ est de
torsion.
\end{defn}
D'après  \cite{Ramachandran} thm 3, $Y$ est une \emph{3-sphère à
homologie entière} si et seulement si $K$ est un corps quadratique
imaginaire dont le groupe de classes est trivial. Alors
$K=\mathbb{Q}(\sqrt{-d})$, où $d$ parcourt l'ensemble
$\{1,2,3,7,11,19,43,67,163\}$.

\begin{prop}\label{revetementcyclique}
Soient $Y$ une 3-sphère à homologie entière et $L/K$ une extension
cyclique de degré $n$. On suppose que $n$ est premier à $2$ dans
tous les cas et premier à $2$ et $3$ pour $K=\mathbb{Q}(\sqrt{-3})$.
Alors on a les suites exactes
$$0\rightarrow \widehat{H}^1(G;Cl(L))\rightarrow \Prod I_{\mathfrak{q}}\rightarrow G\rightarrow0,$$
$$0\rightarrow \widehat{H}^1(G;U_L)\rightarrow \Prod I_{\mathfrak{q}}\rightarrow Cl(L)^G\rightarrow0,$$
où $\Prod I_{\mathfrak{q}}$ désigne le produit des sous-groupes
d'inertie dans $G$ indexés sur l'ensemble des places finies de $K$.
En particulier, si $n=p$ est premier, on a
$$Cl(L)^{C_p}\simeq\mathbb{F}_p^{s-1}.$$
\end{prop}

\begin{proof}
Sous ces hypothèses, le groupe $\widehat{H}^0(G;U_L)$ est trivial.
En effet, $U_K$ est d'ordre~$4$ pour $d=1$, d'ordre~$6$ pour $d=3$
et d'ordre~$2$ sinon. De plus, $\widehat{H}^0(G;U_L)=U_K/N(U_L)$ et
les éléments de ce groupe sont tués par $n$ qui a été choisi premier
au cardinal de $U_K$.

Toutes les différentielles de la suite spectrale $E_r^{**}(X)$ sont
nulles pour $r\geq 2$. Cette suite spectrale est triviale (i.e.
$E_\infty^{**}(X)=E_2^{**}(X)$), convergente et il n'y a que deux
termes non nuls sur chaque diagonale. On obtient ainsi les deux
suites exactes. De plus, pour $n=p$, la première donne
$$s=1+dim_{\mathbb{F}_p}\widehat{H}^1(C_p;Cl(L)).$$
On vérifie la dernière affirmation grâce aux égalités
$$dim_{\mathbb{F}_p}\widehat{H}^1(C_p;Cl(L))=dim_{\mathbb{F}_p}\widehat{H}^0(C_p;Cl(L))=dim_{\mathbb{F}_p}Cl(L)^{C_p},$$
qui proviennent respectivement de la finitude du groupe de classes
$Cl(L)$ et de l'hypothèse $Cl(K)=0$.

\end{proof}

\subsection{Dualité pour la cohomologie équivariante et preuves analogues.}
Nous montrons dans cette section de quelle manière les cohomologies équivariantes ainsi que les suites spectrales
relatives aux $G$-faisceaux $\mathbb{Z}$ et $\mathbb{G}_m$ respectivement, sont liées par une relation de dualité.
Cette dernière provient de dualité d'Artin-Verdier et montre que l'utilisation de  ces deux faisceaux
revient exactement au même. L'utilisation du
$G$-faisceau $\mathbb{Z}$ permet d'obtenir des preuves tout à fait
satisfaisantes du point de vue de la topologie arithmétique mais
l'utilisation du groupe multiplicatif apparait plus naturelle en
cohomologie étale. On conserve les mêmes notations.

\subsubsection{Dualité.}
Ici, le groupe de Galois de l'extension $L/L^G$ est un groupe fini
quelconque, et $L^G$ est totalement imaginaire. On note
$M^D:=Hom(M;\mathbb{Q}/\mathbb{Z})$ le dual d'un groupe abélien de
type fini $M$.

\begin{prop}
On a l'identification $\widehat{H}_G^*(X;\mathbb{Z})\simeq
\widehat{H}_G^*(Z;\mathbb{Z}).$
\end{prop}
\begin{proof}
Le faisceau $\mathbb{Z}$ est adapté (cf \cite{Deninger} 4.6). Le
résultat est vrai d'après la remarque \ref{Z}.

\end{proof}

\begin{lem}\label{lien-1}
On a l'identification
$\widehat{H}_G^n(X;\mathbb{G}_m)\simeq\widehat{H}_G^{n-1}(X;\mathbb{Z}).$
\end{lem}

\begin{proof}
On note toujours $\Omega$ l'ensemble des premiers de $K$ se
ramifiant dans $L$. Pour tout $\mathfrak{q}$ de $\Omega$, on choisit
un premier $\mathfrak{p}$ de $L$ au-dessus de $\mathfrak{q}$. Soient
$\widetilde{Z}_{\mathfrak{p}}$ la composante connexe de
$\widetilde{Z}$ correspondant à $\mathfrak{p}$,
$i:\mathfrak{p}\rightarrow \widetilde{Z}_{\mathfrak{p}}$ l'immersion
fermée, $\eta$ l'inclusion du point générique de
$\widetilde{Z}_{\mathfrak{p}}$ et $G_{\mathfrak{p}}$ le groupe de
décomposition en $\mathfrak{p}$. On a
\begin{equation}\label{lien0}
\widehat{H}_G^n(X;\mathbb{G}_m)\simeq\widehat{H}_G^n(\widetilde{Z};\mathbb{G}_m))
\simeq\Prod_{\Omega}\widehat{H}_{G_{\mathfrak{p}}}^n(\widetilde{Z}_{\mathfrak{p}};\mathbb{G}_m),
\end{equation}
où la deuxième égalité s'obtient en appliquant le lemme de Shapiro
sur le terme initial de la suite spectrale
$\widehat{H}^p(G;H^q(\widetilde{Z};\mathbb{G}_m))\Longrightarrow\widehat{H}^{p+q}_G(\widetilde{Z};\mathbb{G}_m).$
On considère maintenant la suite exacte de
$G_{\mathfrak{p}}$-faisceaux sur $\widetilde{Z}_{\mathfrak{p}}$
\begin{equation}\label{suitelocale}
0\rightarrow\mathbb{G}_{m;\widetilde{Z}_{\mathfrak{p}}}\rightarrow\eta_*\mathbb{G}_{m;\eta}\rightarrow
i_*\mathbb{Z}\rightarrow 0.
\end{equation}
Les faisceaux $R^q(\eta_*)(\mathbb{G}_m)$ sont nuls pour $q\geq1$
(cf \cite{Ma}) et $\eta_*$ préserve les $G_{\mathfrak{p}}$-faisceaux
injectifs. On obtient donc une $G_{\mathfrak{p}}$-résolution
injective de $\eta_*\mathbb{G}_m$ en appliquant $\eta_*$ à une
$G_{\mathfrak{p}}$-résolution injective de $\mathbb{G}_m$. Ceci
permet l'identification
$\widehat{H}_{G_{\mathfrak{p}}}^*(\eta;\mathbb{G}_m)\simeq\widehat{H}_{G_{\mathfrak{p}}}^*(\widetilde{Z}_{\mathfrak{p}};\eta_*\mathbb{G}_m).$
Mais $G_{\mathfrak{p}}$ opère sans inertie sur $\eta$ et
$\mathbb{G}_m$ est adapté sur $\eta$. En effet, $\eta$ est le
spectre de $L_{\mathfrak{p}}$, l'hensélisé de $L$ pour la valuation
donné par $\mathfrak{p}$ (cf \cite{Endler}). Or le groupe de Galois
de ce corps est le même que celui de son complété, qui est de
dimension cohomologie stricte égale à 2. On en déduit
\begin{equation}\label{nul}
\widehat{H}_{G_{\mathfrak{p}}}^*(\eta;\mathbb{G}_m)\simeq\widehat{H}_{G_{\mathfrak{p}}}^*(\widetilde{Z}_{\mathfrak{p}};\eta_*\mathbb{G}_m)=0.
\end{equation}
D'autre part, le foncteur $i_*$ est exact et préserve les
$G_{\mathfrak{p}}$-faisceaux injectifs. On a donc l'identification
\begin{equation}\label{pasnul}
\widehat{H}_{G_{\mathfrak{p}}}^*(\mathfrak{p};\mathbb{Z})\simeq\widehat{H}_{G_{\mathfrak{p}}}^*(\widetilde{Z}_{\mathfrak{p}};i_*\mathbb{Z}).
\end{equation}

La suite exacte longue de cohomologie équivariante associée à
(\ref{suitelocale}) ainsi que les relations (\ref{nul}) et
(\ref{pasnul}) donnent les isomorphismes
\begin{equation}\label{lien1}
\widehat{H}_{G_{\mathfrak{p}}}^{n}(\widetilde{Z}_{\mathfrak{p}};\mathbb{G}_{m;\widetilde{Z}_{\mathfrak{p}}})
\simeq\widehat{H}_{G_{\mathfrak{p}}}^{n-1}(\mathfrak{p};\mathbb{Z}).
\end{equation}

De plus, on a
\begin{equation}\label{lien2}
\Prod_{\Omega}\widehat{H}_{G_{\mathfrak{p}}}^{n-1}(\mathfrak{p};\mathbb{Z})\simeq\widehat{H}_G^{n-1}(Z;\mathbb{Z})\simeq\widehat{H}_G^{n-1}(X;\mathbb{Z}),
\end{equation}
où la première égalité s'obtient à nouveau en appliquant le lemme de
Shapiro sur le terme initial de la suite spectrale (ou directement
sur les groupes de cohomologie équivariante (cf \cite{Br} VII 5).
Finalement, les relations (\ref{lien0}), (\ref{lien1}) et
(\ref{lien2}) permettent de conclure.
\end{proof}

\begin{rem}
On a remarqué dans la preuve précédente que $\mathbb{G}_m$ était
adapté sur $\eta$, où $\eta:=Spec(L_{\mathfrak{p}})$ est le spectre
d'un corps local sur lequel un groupe de Galois $G$ opère. Un tel corps
est de dimension cohomologique $\leq2$ et son groupe de Brauer
s'identifie à $\mathbb{Q}/\mathbb{Z}$ (cf \cite{Milne2}A.1). La
suite spectrale associée ne possède que deux lignes non nulles.
Comme elle converge vers 0, toutes ses différentielles
$$d_3^{p;2}:\widehat{H}^p(G;\mathbb{Q}/\mathbb{Z})\rightarrow
\widehat{H}^{p+3}(G;L_{\mathfrak{p}}^{\times}),$$ sont bijectives.
On retrouve les isomorphismes de la théorie du corps de classe
local.
\end{rem}

\begin{lem}
On a l'identification
$\widehat{H}_G^n(X;\mathbb{Z})\simeq\widehat{H}_G^{n-1}(X;\mathbb{Q}/\mathbb{Z}).$
\end{lem}
\begin{proof}
En considérant la suite spectrale de Leray donnée par l'inclusion du
point générique de $X$, on voit que les groupes $H^q(X;\mathbb{Q})$
sont nuls pour $q\geq1$ (cf \cite{Milne2} II.2.10). Le terme initial
de la suite spectrale est donc nul puisque
$\widehat{H}^*(G;\mathbb{Q})=0$. Ainsi, tous les groupes
$\widehat{H}^q_G(X;\mathbb{Q})$ sont nuls. La suite exacte longue de
cohomologie équivariante associée à la suite exacte de $G$-faisceaux
$$0\rightarrow\mathbb{Z}\rightarrow\mathbb{Q}\rightarrow\mathbb{Q}/\mathbb{Z}\rightarrow0$$ donne le résultat.

\end{proof}

\begin{thm}
On a l'isomorphisme $$\widehat{H}_G^n(X;\mathbb{G}_m)\simeq\widehat{H}_G^{2-n}(X;\mathbb{Z})^D.$$
\end{thm}

\begin{proof}
D'après les deux lemmes précédents, il suffit de montrer l'isomorphisme
\begin{equation} \label{lien-2}
\widehat{H}_G^n(X;\mathbb{Q}/\mathbb{Z})\simeq\widehat{H}_G^{-n}(X;\mathbb{Z})^D.
\end{equation}
De la même manière que dans la preuve du lemme \ref{lien-1},
le théorème de localisation permet de se ramener au spectre d'un corps fini $x=\mathfrak{p}$, sur lequel un groupe de Galois $G$ opère.
Les groupes $H^q(x;\mathbb{Q}/\mathbb{Z})$ sont nuls en dimension supérieure à deux et valent $\mathbb{Q}/\mathbb{Z}$ pour $q=0,1$. La suite spectrale
$$\widehat{H}^p(G;H^q(x;\mathbb{Q}/\mathbb{Z}))\Longrightarrow\widehat{H}_G^p(x;\mathbb{Q}/\mathbb{Z})$$
ne possède que deux ligne non nulles et se réduit (cf \cite{cart} XV 5.11) à la suite exacte longue

$$...\rightarrow\widehat{H}^{n-2}(G;\mathbb{Q}/\mathbb{Z})\rightarrow\widehat{H}^{n}(G;\mathbb{Q}/\mathbb{Z})
\rightarrow\widehat{H}_G^{n}(x;\mathbb{Q}/\mathbb{Z})\rightarrow\widehat{H}^{n-1}(G;\mathbb{Q}/\mathbb{Z})
\rightarrow\widehat{H}^{n+1}(G;\mathbb{Q}/\mathbb{Z})\rightarrow...
$$
D'après \cite{Br} VI.7.3, on a l'identification
$\widehat{H}^{n}(G;\mathbb{Q}/\mathbb{Z})\simeq\widehat{H}^{-1-n}(G;\mathbb{Z})^D.$ En posant $r:=-n$, on obtient

$$...\rightarrow\widehat{H}^{r+1}(G;\mathbb{Z})^D\rightarrow\widehat{H}^{r-1}(G;\mathbb{Z})^D
\rightarrow\widehat{H}_G^{-r}(x;\mathbb{Q}/\mathbb{Z})\rightarrow\widehat{H}^{r}(G;\mathbb{Z})^D
\rightarrow\widehat{H}^{r-2}(G;\mathbb{Z})^D\rightarrow...
$$
De la même manière, en observant la suite spectrale associé au $G$-faisceau $\mathbb{Z}$
sur $x$, on trouve

$$...\rightarrow\widehat{H}^{r-2}(G;\mathbb{Z})\rightarrow\widehat{H}^{r}(G;\mathbb{Z})
\rightarrow\widehat{H}_G^{r}(x;\mathbb{Z})\rightarrow\widehat{H}^{r-1}(G;\mathbb{Z})
\rightarrow\widehat{H}^{r+1}(G;\mathbb{Z})\rightarrow...
$$
Les deux dernières suites exactes se trouvent être duales l'une de
l'autre, ce qui permet d'identifier les groupes
$\widehat{H}_G^n(x;\mathbb{Q}/\mathbb{Z})$ et
$\widehat{H}_G^{-n}(x;\mathbb{Z})^D$. On obtient (\ref{lien-2}) en
appliquant le lemme de Shapiro et le théorème de localisation des
deux cotés, ce qui achève la preuve du théorème.

\end{proof}
D'autre part, les groupes
$H^q(X;\mathbb{Z})$ et $H^{3-q}(X;\mathbb{G}_m)$ sont liés par la
dualité d'Artin-Verdier de la manière suivante (cf \cite{Deninger}):
$$H^3(X;\mathbb{Z})=H^0(X;\mathbb{G}_m)^D=U_L^D,\,\,\, H^2(X;\mathbb{Z})=H^1(X;\mathbb{G}_m)^D=Cl(L)^D,$$
$$H^2(X;\mathbb{G}_m)=H^1(X;\mathbb{Z})^D=0,\,\,\,
H^3(X;\mathbb{G}_m)=H^0(X;\mathbb{Z})^D=\mathbb{Z}^D.$$ De plus,
quel que soit le $G$-module $M$, on a
$\widehat{H}^i(G;M^D)=\widehat{H}^{-1-i}(G;M)^D$ (cf \cite{Br} VI
7.3). On en déduit
$E_2^{p;q}(X;\mathbb{G}_m)=E_2^{-1-p;3-q}(X;\mathbb{Z})^D$.
Vraisemblablement, les différentielles de la première suite
spectrale sont données par les applications duales (ou transposées)
de la deuxième. On obtient alors la proposition suivante par
exactitude du foncteur $Hom(-;\mathbb{Q}/\mathbb{Z})$.

\begin{prop}
Quel que soit $2\leq r\leq\infty$ et quels que soient les entiers $p$ et $q$, on a
$$E_r^{p;q}(X;\mathbb{G}_m)=E_r^{-1-p;3-q}(X;\mathbb{Z})^D$$
\end{prop}

\subsubsection{Analogues topologiques et preuves analogues.}

Soit $M$ une variété topologique compacte de dimension trois,
fermée, lisse, connexe, orientable et sur laquelle le groupe $C_p$
opère fidèlement par automorphismes préservant l'orientation. On
note $Z:=M^{C_p}$ le lieu de ramification du revêtement
$M\rightarrow M/C_p$. Il est constitué de $s$ noeuds, dits noeuds
ramifiés. Soient de plus $H_{tor}(M)$ le sous-groupe de torsion de
$H_1(M;\mathbb{Z})$ et $H_{free}(M)$ le quotient
$H_1(M;\mathbb{Z})/H_{tor}(M)$.

 Respectivement, soit $X$ le spectre de l'anneau d'entiers d'un corps de nombre $L$ sur lequel
$C_p$ opère fidèlement et de sorte que $L^{C_p}$ soit totalement
imaginaire. On note encore $Z$ le lieu de ramification qui est ici
constitué de $s$ places finies ramifiées. Les modules galoisiens
$Cl(L)$ et $U_L/\mu$ sont les analogues arithmétiques de
$H_{tor}(M)$ et $H_{free}(M)$ respectivement. On reprend ci-dessous
le théorème \ref{majoration} ainsi que son analogue topologique.

\begin{thm}\label{majtopar}
En supposant $H_{free}(M/C_p)=0$ dans le cadre topologique, on a les
inégalités suivantes.
$$ s \leq 1+dim_{\mathbb{F}_p}\widehat{H}^1(C_p;H_{free}(M))+
dim_{\mathbb{F}_p}\widehat{H}^0(C_p;H_{tor}(M)).$$
$$ s \leq 1+dim_{\mathbb{F}_p}\widehat{H}^1(C_p;U_L/\mu)+
dim_{\mathbb{F}_p}\widehat{H}^0(C_p;Cl(L)).$$
\end{thm}
On commence par montrer la première inégalité. On réfère à
\cite{Sikora} pour les détails.
\begin{proof}
Le théorème de localisation (cf \cite{Swan} 3.1) fournit
l'isomorphisme
$$\widehat{H}_{C_p}^{n}(M;\mathbb{Z})\simeq
\widehat{H}_{C_p}^{n}(Z;\mathbb{Z}).$$ Le groupe $C_p$ opère
trivialement sur $Z$ donc la suite spectrale $\widehat
H^i(C_p;H^j(Z;\mathbb{Z}))\Longrightarrow\widehat{H}^{i+j}_{C_p}(Z;\mathbb{Z})$
est triviale (cf \cite{Swan} 1.1). Cette suite spectrale ne possède
qu'un seul terme non nul sur chaque diagonale, d'ailleurs isomorphe
à $\mathbb{F}_p^s$. On obtient
$$\widehat{H}_{C_p}^{n}(M;\mathbb{Z})\simeq
\widehat{H}_{C_p}^{n}(Z;\mathbb{Z})\simeq \mathbb{F}_p^s.$$

 En observant le terme initial de la suite spectrale
$\widehat
H^i(C_p;H^j(M;\mathbb{Z}))\Longrightarrow\widehat{H}^{i+j}_{C_p}(M;\mathbb{Z})$,
on obtient immédiatement
$$ s \leq 1+dim_{\mathbb{F}_p}\widehat{H}^0(C_p;H^2(M;\mathbb{Z}))+
dim_{\mathbb{F}_p}\widehat{H}^1(C_p;H^1(M;\mathbb{Z})).$$ La dualité
de Poincaré entraine alors (cf \cite{Sikora} lemme 3.2)
$$ s \leq 1+dim_{\mathbb{F}_p}\widehat{H}^0(C_p;H_1(M;\mathbb{Z}))+
dim_{\mathbb{F}_p}\widehat{H}^1(C_p;H_{free}(M;\mathbb{Z})).$$ On
utilise ensuite l'inégalité
$$dim_{\mathbb{F}_p}\widehat{H}^0(C_p;H_{tor}(M))\leq
dim_{\mathbb{F}_p}\widehat{H}^0(C_p;H_1(M;\mathbb{Z}))$$ qui découle
de l'hypothèse $H_{free}(M/C_p)=0$ (cf \cite{Sikora} preuve du
théorème 1.1(1)).
\end{proof}

La preuve du théorème \ref{majoration}  prend la forme suivante
lorsqu'elle est appliquée au faisceau $\mathbb{Z}$.
\begin{proof}
Le théorème de localisation fournit l' isomorphisme
$$\widehat{H}_{C_p}^{n}(X;\mathbb{Z})\simeq
\widehat{H}_{C_p}^{n}(Z;\mathbb{Z}).$$
 Le groupe $C_p$ opère
trivialement sur $Z$ donc la suite spectrale $\widehat
H^i(C_p;H^j(Z;\mathbb{Z}))\Longrightarrow\widehat{H}^{i+j}_{C_p}(Z;\mathbb{Z})$
est triviale. Cette suite spectrale ne possède qu'un seul terme non
nul sur chaque diagonale, d'ailleurs isomorphe à $\mathbb{F}_p^s$.
On obtient $$\widehat{H}_{C_p}^{n}(X;\mathbb{Z})\simeq
\widehat{H}_{C_p}^{n}(Z;\mathbb{Z})\simeq \mathbb{F}_p^s.$$

 En observant le terme initial de la suite spectrale
$\widehat
H^i(C_p;H^j(X;\mathbb{Z}))\Longrightarrow\widehat{H}^{i+j}_{C_p}(X;\mathbb{Z})$,
on obtient immédiatement
$$ s \leq dim_{\mathbb{F}_p}\widehat{H}^1(C_p;H^2(X;\mathbb{Z}))+
dim_{\mathbb{F}_p}\widehat{H}^0(C_p;H^3(X;\mathbb{Z})).$$ La dualité
d'Artin-Verdier entraine alors
$$ s \leq dim_{\mathbb{F}_p}\widehat{H}^0(C_p;Cl(L))+
dim_{\mathbb{F}_p}\widehat{H}^1(C_p;U_L).$$ L'inégalité
\ref{deuxmajoration}
$$dim_{\mathbb{F}_p}\widehat{H}^0(C_p;U_L)\leq 1+
dim_{\mathbb{F}_p}\widehat{H}^0(C_p;U_L/\mu)$$ permet de conclure.

\end{proof}

On reprend maintenant le théorème \ref{minoration} ainsi que son
analogue topologique.
\begin{thm}\label{mintopar}
Si $H_1(M;\mathbb{Z})$ s'identifie à $H_{free}(M)\oplus H_{tor}(M)$
en tant que $C_p$-module et si $s\geq1$, alors $$s\geq
1+dim_{\mathbb{F}_p}\widehat{H}^0(C_p;H_{tor}(M)).$$ Si $Cl(K)=0$,
alors $$s\geq 1+dim_{\mathbb{F}_p}\widehat{H}^0(C_p;Cl(L)).$$
\end{thm}

De la même manière, la démonstration du théorème \ref{minoration}
prend une forme plus agréable lorsqu'elle est appliquée au faisceau
$\mathbb{Z}$. Les preuves de ces deux résultats reposent alors sur
le théorème de localisation, la dualité de Poincaré (respectivement
d'Artin-Verdier) et sur une étude un peu plus fine de la suite
spectrale (section 3.2 de \cite{Sikora} et section 4.2 de ce papier)
permettant de montrer que certains modules du terme initial
survivent à l'infini.

\subsection{Cas des corps de nombres admettant des plongements réels.}
Dans cette section, les corps de nombres $L$ et $K=L^G$ ne sont plus
nécessairement totalement imaginaires. Il est alors indispensable
d'utiliser la topologie étale d'Artin-Verdier qui tient compte des
places à l'infini. Afin d'alléger l'exposé, nous présentons ici
quelques résultats sans démonstrations pour lesquelles nous
renvoyons à \cite{moi}.

Soient $X$ le spectre de l'anneau d'entiers de $L$, $X_{\infty}$
l'ensemble des places archimédiennes de $L$, $\overline{X}$ le
couple $(X;X_{\infty})$ et $j$ l'inclusion de $X$ dans
$\overline{X}$. On munit $\overline{X}$ de la topologie étale
d'Artin-Verdier (cf \cite{Bienenfeld}). Les groupes de cohomologie
du groupe multiplicatif
$\mathbb{G}_{m;\overline{X}}:=j_*\mathbb{G}_m$ sont respectivement
$U_L$, $Cl(L)$, 0, $\mathbb{Q}/\mathbb{Z}$, et 0 en dimension
supérieure à trois (cf \cite{Bienenfeld} 2.7). Le groupe de Galois
de l'extension $L/K$ opère sur $\overline{X}$ et l'on définit la
catégorie des $G$-faisceaux pour cette topologie ainsi que les
groupes de cohomologie équivariante
$\widehat{H}^*_G(\overline{X};F)$. Les résultats des sections 2 et 3
peuvent être appliqués en prenant certaines précautions.

Si $\mathfrak{p}$ est une place finie de $L$, les notations
intervenant dans la proposition suivante sont celles qui ont été
introduites dans la section 4.1. Si $\mathfrak{p}$ est une place
archimédienne complexe, le couple
$(\overline{L};\overline{\mathfrak{p}})$ est un hensélisé du corps
valué $(L;\mathfrak{p})$ (cf \cite{Endler}), où
$\overline{\mathfrak{p}}$ désigne un prolongement de $\mathfrak{p}$
à $\overline{L}$. On pose alors
$U_{\mathfrak{p}}:=\overline{L}^{\times}$ et
$I_{\mathfrak{p}}^{(ab)}:=I_{\mathfrak{p}}=G_{\mathfrak{p}}\subseteq
G_{\overline{\mathfrak{p}}}\simeq \mathbb{Z}/2\mathbb{Z}$.

\begin{prop} \label{calculs sur overline}
Soit $G$ un groupe fini opérant fidèlement sur $\overline{X}$. On a
les isomorphismes ci-dessous, où les deux premiers produits sont
pris sur toutes les places de $K$ et le troisième sur l'ensemble des
places finies de $K$.

\begin{itemize}
\item $\widehat{H}_G^*(\overline{X};\mathbb{G}_{m;\overline{X}})\simeq\Prod\widehat{H}^*(G_{\mathfrak{p}};U_{\mathfrak{p}}).$
\item  $\widehat{H}_G^0(\overline{X};\mathbb{G}_{m;\overline{X}})\simeq\Prod I_{\mathfrak{q}}^{(ab)}.$
\item  $\widehat{H}_G^1(\overline{X};\mathbb{G}_{m;\overline{X}})\simeq\Prod\mathbb{Z}/e_{\mathfrak{q}}\mathbb{Z}.$
\end{itemize}
\end{prop}

Tous les résultats de la section 4.3 se généralisent ainsi à toutes
les extensions de corps de nombres avec de légères modifications.
Par exemple, si $L/K$ est une extension de corps de nombres de
groupe $C_p$, on a la proposition suivante, où $s$ (respectivement
$\overline{s}$) désigne le nombre de places finies (respectivement
finies et infinies) ramifiées.
\begin{prop}\label{majmininfini}
On a la majoration $ s \leq
1+dim_{\mathbb{F}_p}\widehat{H}^1(C_p;U_L/\mu)+
dim_{\mathbb{F}_p}\widehat{H}^0(C_p;Cl(L)).$

Si le groupe $Cl(K)$ est trivial, alors $\overline{s}\geq
1+dim_{\mathbb{F}_p}\widehat{H}^0(C_p;Cl(L)).$
\end{prop}

On considère à nouveau l'action (fidèle) du groupe $C_p$ sur une
3-variété $M$ et sur un corps de nombres $L$ respectivement, en
gardant les notations (et les hypothèses) précédentes. On note aussi
$s_0$ le nombre de places finies ramifiées dans cette extension.
\begin{prop}\label{egaltopar}
Si $H_{free}(M)=0$ et si $s\geq1$, alors
$$H_{tor}(M)^{C_p}\simeq\mathbb{F}_p^{s-1}.$$

Si $U_L/\mu=0$, alors $p=2$, $s_0\geq1$ et
$$Cl(L)^{C_2}\simeq\mathbb{F}_2^{s_0 -1}.$$
\end{prop}

\begin{proof}
Dans le cadre topologique, les hypothèses $H_{free}(M)=0$ et
$s\geq1$ assurent que la suite spectrale $E_2^{**}(M;\mathbb{Z})$
est triviale (cf \cite{Sikora} 3.3 et 3.4). De la même manière,
l'hypothèse $U_L/\mu=0$ permet de montrer que la suite spectrale
$E_2^{**}(\overline{X};\mathbb{G}_m)$ est triviale. Il faut ici
utiliser les résultats de la section 4.2 généralisés aux corps de
nombres quelconques et considérer de plus le morphisme de suites
spectrales induit par l'inclusion du point générique de
$\overline{X}$ (l'ensemble des valuations de $L$), pour montrer que
toutes les différentielles sont nulles (cf \cite{moi}).

Le fait que les deux suites spectrales soient triviales donne
immédiatement les deux résultats précédents.
\end{proof}


\subsection{Complément technique.}
\begin{lem} \label {clef}
Soient $L/K$ une extension galoisienne de corps de nombres de groupe
$G$, $D$ et $R$ les anneaux d'entiers de $L$ et $K$ respectivement.
On suppose que $Cl(K)$ est trivial. On note $U_L$ (respectivement
$U_K$) le groupe des unités de $L$ (respectivement de $K$). Soient
$(\mathfrak{q}_j)_{1\leq j\leq r}$ les idéaux premiers finis de $K$
qui se ramifient dans l'extension $L/K$ et $(\mathfrak{p}_i)_{1\leq
i\leq s}$ les idéaux premiers de $L$ au-dessus des $\mathfrak{q}_j$.
Soit $U_{\mathfrak{p}_i}$ (respectivement $U_{\mathfrak{q}_j}$) le
groupe des unités de l'hensélisé $D_{\mathfrak{p}_i}^h$ de l'anneau
local $D_{(\mathfrak{p}_i)}$ (respectivement de l'hensélisé
$R_{\mathfrak{q}_j}^h$ de l'anneau local $R_{(\mathfrak{q}_j)}$),
$L_i$ et $K_j$ les corps de fractions de $D_{\mathfrak{p}_i}^h$ et
$R_{\mathfrak{q}_j}^h$.

Si $Cl(K)$ est trivial, alors le morphisme canonique
$H^1(G;U_L)\rightarrow H^1(G;\prod_{1\leq i\leq
s}U_{\mathfrak{p}_i})$ est injectif.
\end{lem}
\begin{proof}
Soit $X^0$ l'ensemble des places finies de $L$. Le groupe de Galois
opère sur $X^0$ et on a le morphisme de $\mathbb{Z}[G]$-modules
$div:L^{\times}\rightarrow \sum_{\mathfrak{p}\in X^0}\mathbb{Z}$
donné par les valuations. On note $W$ l'image du morphisme $div$. On
a le morphisme de suites exactes de $\mathbb{Z}[G]$-modules:
\[ \xymatrix{
0\ar[d]\ar[r] &U_L\ar[d]\ar[r]&L^{\times}  \ar[d]\ar[r]^{div}&W \ar[d]^{p}\ar[r]&0\ar[d]\\
0 \ar[r]      &\prod_{{1\leq i\leq s}} U_{\mathfrak{p_i}}\ar[r]           &\prod_{1\leq i\leq s} L^{\times}_i  \ar[r]      &\prod_{1\leq i\leq s} \mathbb{Z} \ar[r]      &0\\
} \] En utilisant le théorème de Hilbert 90 et le fait que les
$H^q(G;...)$ forment un foncteur cohomologique, on obtient le
morphisme de suites exactes suivant.
\[ \xymatrix{
0\ar[d]\ar[r] &U_K\ar[d]\ar[r]&K^{\times}  \ar[d]\ar[r]^{div\mid K^{\times}}&W^G \ar[d]^{p\mid W^G}\ar[r]&H^1(G;U_L)\ar[d]^j\ar[r]&0\ar[d]\\
0 \ar[r]   &\prod_{{1\leq j\leq r}} U_{\mathfrak{q_j}}\ar[r]     &\prod_{1\leq j\leq r} K^{\times}_j  \ar[r]^{\overline{v}}      &\prod_{1\leq j\leq r} \mathbb{Z} \ar[r]   & H^1(G;\prod_{1\leq i\leq s}U_{\mathfrak{p}_i}) \ar[r]     &0\\
} \] Alors $j$ s'identifie au morphisme $$W^G/Im(div\mid
K^{\times})\longrightarrow (\prod_{1\leq i\leq r}
\mathbb{Z})/Im(\overline{v})$$ induit par $p$, dont on montre
facilement qu'il est injectif lorsque $Cl(K)$ est trivial.

\end{proof}

\section{Conclusion.}
Nous donnons ici une interprétation des résultats et de leurs
preuves exposés dans ce travail, afin d'essayer d'éclaircir et
d'approfondir le dictionnaire de la topologie arithmétique.

\subsection{Quelques éléments du dictionnaire.}

\subsubsection{}

Dans les versions précédentes du dictionnaire de la topologie
arithmétique (\cite{Ramachandran} et \cite{Reznikovpreprint}), le
groupe de Galois $G_L^{nr}$ de l'extension maximale non ramifiée en
toutes les places (archimédiennes et ultramétriques) d'un corps de
nombre $L$ est vu comme l'analogue du groupe fondamental topologique
de la variété "correspondante" $M$. En effet, il s'agit du groupe
fondamental pour la topologie étale d'Artin-Verdier. En suivant
cette idée, $H_1(M;\mathbb{Z})=\pi_1(M)^{ab}$ devrait être
l'analogue du groupe de Galois de l'extension abélienne maximale non
ramifiée de $L$. Par la théorie du corps de classe, ce groupe
s'identifie à $Cl(L)$. Il n'y aurait donc pas de place pour le
groupe des unités. Cette même analogie conduit à montrer que la
version arithmétique de la conjecture de Poincaré est fausse (cf
\cite{Ramachandran}).

C'est peut-être pour éviter des contradictions de ce type que A.
Reznikov a proposé dans la deuxième version du dictionnaire
\cite{Reznikov} 12, de voir un corps de nombre $L$ comme une
3-variété $M$ bordant une 4-variété $N$ de sorte que le morphisme
$\pi_1(M)\rightarrow\pi_1(N)$ soit surjectif, pour considérer les
invariants cohomologiques de $N$ (et non de $M$). Cependant, le
non-sens évoqué ci-dessus est toujours présent dans ce deuxième
dictionnaire.


Comme nous le précisons ci-dessous, le travail exposé dans ce papier
confirme clairement la première version \cite{Reznikovpreprint} du
dictionnaire, et écarte la deuxième \cite{Reznikov}. En effet,
toutes les preuves sont basées sur la cohomologie de $\overline{X}$
et de $M$ pour aboutir à des résultats respectant parfaitement le
dictionnaire \cite{Reznikovpreprint}.

\subsubsection{}

Nous soutenons les analogies suivantes.

Les corps de nombres, les extensions galoisiennes de corps de
nombres $L/L^G$, les places finies, les places finies ramifiée et
l'ensemble des places finies ramifiés doivent être vus
respectivement comme des 3-variétés, des revêtements ramifiés
galoisiens de 3-variétés $M\rightarrow M/G$, des noeuds dans $M$,
des noeuds ramifiés et comme le lieu de ramification (plus
précisément, les sous-groupes de décomposition et d'inertie se
correspondent). Ici, le terme \emph{3-variété} doit probablement
prendre un sens plus large.

De plus, $Cl(L)$ et $U_L/\mu$ correspondent, de manière compatible à
l'action d'un groupe de Galois, aux groupes $H_{tor}(M)$ et
$H_{free}(M)$. D'ailleurs, le paragraphe suivant confirme à nouveau
ces deux dernières correspondances.

\subsection{Cohomologie de S. Lichtenbaum et topologie arithmétique.}

\subsubsection{}

Considérons la topologie Weil-étale. Les résultats de S. Lichtenbaum
dans \cite{Lichtenbaum} montrent que les groupes de cohomologie à
support compact de $\overline{X}$ à coefficients entiers devraient
être les suivants (cf \cite{Lichtenbaum} 6.3).
$$
\begin{array}{rcl}
H_c^{q}(\overline{X};\mathbb{Z})&=& 0\,\,\,\, pour\,\, q=0,\\
          &=& (\prod_{X_{\infty}}\mathbb{Z})/\mathbb{Z}\,\,\,\, pour\,\, q=1,\\
&=&Pic^1(\overline{X})^{\mathcal{D}}\,\,\,\, pour\,\, q=2,\\
&=&(\mu_L)^{\mathcal{D}}\,\,\,\, pour\,\, q=3,\\
&=&0\,\,\,\, pour\,\, q\geq4.\\
\end{array}$$
Ci-dessus, $Pic^1(\overline{X})$ est le noyau du morphisme du groupe
de classe d'Arakelov dans $\mathbb{R}^{\times}$ donné par la valeur
absolue. On note aussi
$A^{\mathcal{D}}:=Hom(A;\mathbb{R}/\mathbb{Z})$ le dual de
Pontryagin d'un groupe topologique abélien séparé et localement
compact $A$. D'autre part, le sous-groupe de torsion et le quotient
libre de $H^2_c(\overline{X};\mathbb{Z})$ sont donnés par la suite
exacte (cf \cite{Lichtenbaum} 6.4)
$$0\rightarrow Cl(L)^{\mathcal{D}}\rightarrow Pic^1(\overline{X})^{\mathcal{D}}\rightarrow Hom(U_L;\mathbb{Z})\rightarrow0.$$

D'après l'interprétation de C. Deninger dans (\cite{Deninger3} 7),
la cohomologie de Lichtenbaum correspond à la cohomologie des
faisceaux sur $M$, qui s'identifie d'ailleurs à la cohomologie
singulière, car une variété topologique est un espace localement
contractile. Puisque $M$ est une variété compacte, la cohomologie à
support compact s'identifie à la cohomologie usuelle. Ainsi,
$H^2(M;\mathbb{Z})$ doit être l'analogue du groupe
$Pic^1(\overline{X})^{\mathcal{D}}$. Lorsque $M$ est orientable, la
dualité de Poincaré identifie $H^2(M;\mathbb{Z})$ et
$H_1(M;\mathbb{Z})$ en tant que groupes abéliens. Dans la situation
équivariante, l'action du groupe de Galois sur ces deux groupes se
trouve être inversée à travers cet isomorphisme. De plus,
$H_c^{1}(\overline{X};\mathbb{Z})$ s'identifie à
$Hom(U_L/\mu;\mathbb{Z})$, c'est-à-dire l'analogue de
$Hom(H_{free}(M;\mathbb{Z});\mathbb{Z})=Hom(\pi_1(M);\mathbb{Z})$.

\subsubsection{D'autres éléments du dictionnaire.}

D'après le paragraphe précédent, les analogues arithmétiques des
groupes $H^1(M;\mathbb{Z})$ et $H^2(M;\mathbb{Z})$ devraient être
$Hom(U_L/\mu;\mathbb{Z})$ et $Pic^1(\overline{X})^{\mathcal{D}}$
respectivement. Si l'on voit $\overline{X}$ comme un espace
orientable, le groupe $H_1(M;\mathbb{Z})$ devrait être l'analogue de
$Pic^1(\overline{X})^{\mathcal{D}}$. L'action d'un éventuel groupe
de Galois $G$ sur ce dernier devrait se faire à travers le groupe
opposé à $G$, afin d'être compatible à l'action de $G$ sur
$H_1(M;\mathbb{Z})$. De plus, $Cl(L)$ et $U_L$ correspondent
respectivement aux groupes $H_{tor}(M)$ et $H_{free}(M)$ déduits de
$H_1(M;\mathbb{Z})$. Une sphère à homologie entière (respectivement
rationnelle) est un corps de nombres dont l'analogue de
$H_1(M;\mathbb{Z})$ est nul (respectivement de torsion).

Enfin, le groupe $G_L^{nr}$ ne serait pas l'analogue de $\pi_1(M)$
(puisque les abélianisés de ces groupes ne se correspondent pas),
mais nous espérons pouvoir revenir plus tard sur cette question.
\begin{rem}
Soit $\overline{X}\rightarrow\overline{Y}$ un morphisme au-dessus de
$\overline{Spec(\mathbb{Z})}$ donné par une extension de corps de
nombre $L/K/\mathbb{Q}$. L'inclusion $I_K\hookrightarrow I_L$ induit
morphisme $Pic^1(\overline{X})^{\mathcal{D}}\rightarrow
Pic^1(\overline{Y})^{\mathcal{D}}$. Ainsi, ce groupe dépend
fonctoriellement de $\overline{X}$ et de manière covariante.
\end{rem}

\subsection{Comparaison des hypothèses.}

Les résultats de A. Sikora étudiés dans ce travail donnent une
majoration (\ref{majtopar}), une minoration (\ref{mintopar}) et une
égalité dans un cas particulier (\ref{egaltopar}), du nombre $s$ de
places (respectivement de noeuds) ramifiées dans un revêtement
cyclique d'ordre premier. Pour la majoration et l'égalité, les
hypothèses faites en topologie sont strictement plus fortes qu'en
arithmétique. De plus, d'après (\cite{Sikora} 5), elles sont toutes
nécessaires en topologie. Malheureusement, nous n'avons pas été en
mesure de démontrer (à l'aide des mêmes méthodes) la minoration de
$s$ en arithmétique à partir de l'analogue des hypothèses
topologiques (qui peuvent être formulées grâce au paragraphe
précédent).

Néanmoins, ce travail donne le sentiment que le cadre topologique
est "plus général" que celui des corps de nombres. En effet,
l'arithmétique apparaît ici beaucoup plus rigide que le cadre
topologique. Cette idée vague est d'ailleurs nettement confirmée
dans \cite{Deninger3}. On peut aussi illustrer ce fait par les
observations suivantes. Un corps de nombres est de manière unique un
revêtement de $\mathbb{Q}$, alors qu'une 3-variété orientable peut
s'obtenir d'un grand nombre de manières différentes comme revêtement
de $\mathbb{S}^3$. Il existe exactement dix sphères à homologie
entière en arithmétique (toutes de degré un ou deux) et une infinité
en topologie (cf \cite{Ramachandran}). L'analogue de
$H_1(M;\mathbb{Z})$ est plus gros que le groupe de Galois de
l'extension abélienne maximale non ramifiée de $L$.

Il est d'ailleurs amusant d'imaginer à quoi ressembleraient les
preuves des mêmes résultats basées sur la topologie Weil-étale. En
supposant que la ligne
$E_2^{i;3}(\overline{X};\varphi_!\mathbb{Z})=\widehat{H}^i(C_p;\mu_L)$
soit non nulle et qu'elle survive à l'infini, on obtiendrait
exactement les mêmes démonstrations dans les deux contextes.

On se rend alors compte que les différences entre les preuves
arithmétiques et topologiques que nous avons proposées proviennent
des "défauts" de la cohomologie étale. Cependant, lorsqu'il n'y pas
de ramification à l'infini, ces mêmes défauts n'apparaissent pas
dans les groupes de cohomologie étale équivariante modifiée. En
effet, ces derniers donnent les "bons" groupes liés à la
ramification (c'est à dire les mêmes que dans le cadre topologique),
ce qui renforce l'analogie des preuves que nous avons proposées. On
vérifie cette affirmation en utilisant le théorème de localisation
et en observant que les termes initiaux des suites spectrales (par
exemple à coefficients dans $\mathbb{Z}$), définies sur le spectre
d'un corps fini et sur un cercle munis d'une action d'un groupe fini
$G$, sont en fait les mêmes. Cette analogie n'est d'ailleurs pas
respectée par les (quatre premiers) groupes de cohomologie étale
équivariante non modifiée. En effet, si $G$ opère trivialement, on a
$H^1(Spec(\mathbb{F}_q)_{et};G;\mathbb{Z})=0$ et
$H^1(\mathbb{S}^1;G;\mathbb{Z})=\mathbb{Z}$, comme le montre la
suite spectrale \cite{tohoku} 5.2.9.

Lorsqu'il y a de la ramification à l'infini, les places complexes
apparaissent comme des points dans la cohomologie équivariante
modifiée.

\subsection{Places finies et places archimédiennes.}

\subsubsection{}
Dans le topos associé à la topologie Weil-étale sur $\overline{X}$,
une place finie est donnée par l'inclusion fermée du topos
classifiant $B_{\mathbb{Z}}$. Par ailleurs, l'espace classifiant du
groupe discret $\mathbb{Z}$ est le cercle $\mathbb{S}^1$, dont le
topos est homotopiquement équivalent à $B_{\mathbb{Z}}$ (cf
\cite{Moerdijk} IV.1.1). Il est donc naturel de voir géométriquement
une place finie comme l'immersion fermée d'un cercle.
\subsubsection{}
D'après N. Ramachandran, on sait que $\overline{X}$ doit être vu
comme la compactification d'une variété non compacte correspondant à
$X$.

Dans le topos associé à la topologie étale (d'Artin-Verdier) sur
$\overline{X}$, une place archimédienne est donnée par l'inclusion
du topos ponctuel, c'est-à-dire un point du topos étale (cf
\cite{SGA4} IV.6). C'est la raison pour laquelle les calculs de
cohomologie étale font apparaître les places archimédiennes comme
des points, les bouts d'une 3-variété non compacte.

Cette analogie offre quelques contradictions. D'une part, toute
variété orientable de dimension trois se réalise comme revêtement
ramifié de $\mathbb{S}^3$ dont le lieu de ramification est composé
de noeuds. Or un corps de nombres non totalement réel est, de
manière unique, un revêtement de $\mathbb{Q}$ ramifié à l'infini.
Par exemple, une extension quadratique imaginaire de $\mathbb{Q}$
ramifiée en au moins un premier ne peut pas être vue comme une
variété orientable $M$ au-dessus de $M/C_2=\mathbb{S}^3$. En effet,
dans cette situation, le lieu de ramification ne peut pas être
constitué de noeuds et d'un point (isolé dans le lieu de
ramification). D'autre part, les places archimédiennes doivent être
comptées comme composantes du lieu de ramification dans
\ref{majmininfini}, et doivent être ignorées dans \ref{egaltopar}.

Il semble donc que le fait de voir un corps de nombres comme une
simple variété orientable pose quelques problèmes dans
l'interprétation des places archimédiennes. D'ailleurs, dans les
travaux de C. Deninger (cf \cite{Deninger2} et \cite{Deninger3}), il
n'est pas clair qu'un corps de nombres corresponde à une variété
proprement dite.

De plus, dans le topos associé à la topologie Weil-étale, une place
archimédienne devrait être donnée par l'inclusion du topos
classifiant $B_{\mathbb{R}}$ du groupe topologique $\mathbb{R}$, qui
est loin d'être équivalent au topos ponctuel.

\subsubsection{}
Toujours dans les travaux de Christopher Deninger, une place
archimédienne correspond à un point sur lequel $\mathbb{R}$ opère
trivialement, ce qui donne lieu au topos $B_{\mathbb{R}}$ (cf
\cite{SGA4} IV. 2.5). Respectivement, une place finie $\mathfrak{p}$
correspond à l'immersion fermée d'un cercle
$\mathbb{R}/log(N\mathfrak{p})\mathbb{Z}=\mathbb{R}/l(\mathfrak{p})\mathbb{Z}$
sur lequel $\mathbb{R}$ opère naturellement. Le topos associé est le
topos induit sur l'objet $\mathbb{R}/l(\mathfrak{p})\mathbb{Z}$ de
$B_{\mathbb{R}}$ (cf \cite{SGA4} IV.5.1). Ce dernier topos est
canoniquement équivalent à $B_{l(\mathfrak{p})\mathbb{Z}}$ (cf
\cite{SGA4} IV.5.8). Par ailleurs, dans le topos associé à la
topologie Weil-étale, une place finie $\mathfrak{p}$ est donnée plus
exactement par l'immersion fermée de $B_{W_{k(\mathfrak{p})}}$. De
plus, en normalisant convenablement, le logarithme de la valeur
absolue donne
$W_{k(\mathfrak{p})}=log(N\mathfrak{p})\mathbb{Z}=l(\mathfrak{p})\mathbb{Z}\subseteq\mathbb{R}$.
Ces deux points de vue sont donc compatibles.


\end{document}